\documentclass[10pt,a4paper]{article}
\usepackage{lineno,hyperref}
\modulolinenumbers[5]
\usepackage[utf8]{inputenc}
\usepackage{amsmath}
\usepackage{amssymb}
\usepackage{graphicx}
\usepackage{subfig}
\usepackage{xcolor,soul}
\usepackage[top=2cm,bottom=2cm,left=1.5cm,right=1.5cm]{geometry}

\newtheorem{remark}{Remark}

\date{}

\begin{document}
\title{A fully non-invasive hybrid IGA/FEM scheme for the analysis of localized non-linear phenomena}

\maketitle

\begin{center}
	Evgeniia Lapina$^{1,2}$, Paul Oumaziz$^1$, Robin Bouclier$^{1,2}$, Jean-Charles Passieux$^1$ \\
	~\\
	$^1$ Institut Clement Ader (ICA), Université de Toulouse, CNRS-INSA-UPS-ISAE-Mines Albi, 5 rue Caroline Aigle, Toulouse, 31400, France \\
	$^2$ Institut de Mathématiques de Toulouse (IMT), Université de Toulouse, CNRS-INSA-UT1-UT2-UPS,118, route de Narbonne , Toulouse,F-31062 Cedex 9, France\\
    ~\\
	evgeniia.lapina@insa-toulouse.fr, paul.oumaziz@insa-toulouse.fr, bouclier@insa-toulouse.fr,passieux@insa-toulouse.fr
\end{center}

\textbf{Abstract:} This work undertakes to combine the interests of IsoGeometric Analysis (IGA) and standard Finite Element Methods (FEM) for the global/local simulation of structures. The idea is to adopt a hybrid global-IGA/local-FEM modeling, thereby benefiting from: (i) the superior geometric description and per-Degree-Of-Freedom accuracy of IGA for capturing global, regular responses, and (ii) the ability of FEM to compute local, strongly non-linear or even singular behaviors. For the sake of minimizing the implementation effort, we develop a coupling scheme that is fully non-invasive in the sense that the initial global spline model to be enriched is never modified and the construction of the coupling operators can be performed using conventional FE packages. The key ingredient is to express the FEM-to-IGA bridge, based on Bézier extraction, to transform the initial global spline interface into a FE one on which the local FE mesh can be constructed. This allows to resort to classic FE trace operators to implement the coupling. It results in a strategy that offers the opportunity to simply couple an isogeometric code with any robust FE code suitable for the modelling of complex local behaviors. The method also easily extends in case the users only have at their disposal FE codes. This is the situation that is considered for the numerical illustrations. More precisely, we only make use of the FE industrial software Code\_Aster to perform efficiently and accurately the hybrid global-IGA/local-FEM simulation of structures subjected locally to cracks, contact, friction and delamination.

\medskip

\textbf{Keywords:} Isogeometric Analysis, Non-invasive global/local coupling, Multiscale, Bézier extraction, Industrial FE code

\section{Introduction}

Originally introduced in~\cite{hughes05,cottrell09}, the core idea of IsoGeometric Analysis (IGA) is to resort to the same higher-order and smooth bases, in particular made of B-Spline and Non-Uniform-Rational-B-Spline (NURBS) functions~\cite{Cohen80,Piegl97}, for the representation of the geometry in Computer-Aided Design (CAD) as well as for the approximation of solutions fields in numerical simulations. All over the paper, we employ the terminologies spline or IsoGeometric (IG) indifferently to denote a NURBS or a B-Spline object. The use of such functions quickly made IGA highly attractive with respect to the standard Finite Element Method (FEM) for two main reasons. On the one hand, a common geometrical model can be used by both the designers and analysts, thereby facilitating the dialog between their two worlds. On the other hand, the (possibly more) regular approximation spaces offered by IGA may be included into the $C^0$ spaces given by classic FEM (provided that both IGA and FEM come with the same polynomial degree). As a result, IGA can be interpreted as a projection of FEM onto a regular reduced basis~\cite{Tirvaudey19}. This interpretation in terms of reduced FEM modeling accounts for the increased per-Degree-Of-Freedom (per-DOF) accuracy attributed to IGA: in some sense, IGA allows to capture a regular solution as accurately as FEM but with fewer DOF because its discretization space is included into the FE $C^0$ space. Hence, this technology is now often seen as a high-performance computational tool in the community.

However, the latter point of view also highlights that IGA, in its standard brute form, does not appear suitable when the solution is not regular anymore. For instance, this often occurs at the local scale in structural mechanics (\emph{e.g.}, cracks, local contact, delamination, or local heterogeneities which involves displacement and/or strain discontinuities). As a result, numerous sophisticated methods have been developed over the years to make possible local simulations with IGA. For the representation of geometrical details, one solution, to avoid tedious (or even impossible) spline re-parametrizations leading to the splitting of the geometry into several $C^0$ patches~\cite{Xu13,Massarwi19}, may be to resort to immersed IGA where the geometry is given in terms of trimming entities while the numerical approximation space is built on an embedding spline cuboid~\cite{Ruess14,Wei21,Wang21}. Then, regarding for instance fracture and/or delamination, one may refer to the IG version of XFEM, namely XIGA~\cite{Luycker11,Yuan21,Fathi21}, or to IG cohesive elements~\cite{Verhoosel11,Dimitri14}, or even to phase-field approaches~\cite{Borden14,Proserpio20,Paul20}, to name a few. Eventually, all these methods seem to have a very high level of complexity and therefore may require significant effort to be understood, and implemented from a classic IG code. On the contrary, FEM appears adapted to simulate local, strongly non-linear or even singular behaviors due to its reduced regularity and its meshing freedom. Furthermore, FEM benefits from more that 50 years of developments and practices so numerous enhanced FE implementations, both efficient and robust, exist to simulate various local behaviors.

In this context, we propose here to adopt a hybrid global-IGA/local-FEM modeling, so as to end up with a combined strategy which mixes the interests of both analysis technologies for global/local simulations: efficiency of IGA for geometric description and for capturing global, regular response; and, ability of FEM to compute local, strongly non-linear or even singular behaviors. In addition, we seek for the simplest possible strategy in terms of implementation effort in order to be able to use any (possibly industrial) FE code for the local behavior. Consequently, our starting point is to consider what is now referenced as the non-invasive global/local coupling methodology in the field. In an iterative process, a part of the global model is replaced by the more detailed local model exactly and non-invasively: the global model is never modified; only interface displacements and reaction forces are exchanged. This strategy along with its concept of non-invasiveness have been successfully applied in FEM and are still gathering a considerable interest in the community (see~\cite{Gendre09} for local plasticity, ~\cite{Passieux13,Duarte21,Oumaziz_chilien22,Meray22} for crack propagation,~\cite{Gerasimov_18,Noii_20,Aldakheel_21} for fracture modeling with the phase-field approach, \cite{Duval16,Gosselet18} for domain decomposition solvers, \cite{Oumaziz19} for multi-contact problems,~\cite{Guinard18} for real aeronautical structures, and \cite{Wangermez20} for multiscale periodic heterogeneous materials, to name a few). In this work where we consider the coupling of a global IG model with a local FE model, the non-invasive global/local framework appears even more relevant~\cite{Bouclier16}: (i) it naturally avoids costly spline re-parametrization procedures, which may have been necessary otherwise to incorporate a truly-independent local region within the initial IG model, (ii) the global IG stiffness operator can be assembled and factorized only once and the IG system to be solved remains well-conditioned regardless of the shape of the local region, and (iii) it offers the opportunity to simply couple an IG code with any existing robust FE code suitable for the modelling of complex local behaviors.

The remaining difficulty when considering domain coupling within IGA is the formulation and implementation of a possibly non-conforming coupling. Inspired from immersed methods~\cite{Ruess14,Wei21,Wang21}, the usual coupling of the non-invasive strategy by means of Lagrange multipliers was replaced by a Nitsche-based coupling to answer this issue in the field of full global/local IGA~\cite{Bouclier18,Bouclier_book_22}. However, such an approach appears inconsistent with the use of standard industrial FE codes. As a remedy, we propose here to make use of the FE-type description of the local model: the idea is to call upon efficient (classic) FE meshing procedures to recover the simple case of a conforming interface, \emph{i.e.} to align the boundary of the local FE mesh to the edges of some global (knot-span) elements. Then, the second ingredient consists in resorting to the existing FEM-to-IGA bridge~\cite{Tirvaudey19,Colantonio20}, based on Bézier or Lagrange extraction operators~\cite{Borden11,Scott11,Schillinger16}, to transform the initial interface within the global IG model into a FE interface on which the local FE mesh can be constructed. With all these elements, the actual coupling between IGA and FEM can be done explicitly by resorting to only standard FE trace operators. These coupling operators being most times available in FE codes, the coupling can be carried out using only FE industrial packages. In the end, we arrive at a fully non-invasive strategy in the sense that not only the global/local coupling is non-invasive but also the construction of the coupling operators. Although the theory applies for any higher-order B-splines and NURBS, we will restrict ourselves to the quadratic case for the numerical experiments since we seek to use only (industrial) FE packages for the implementation, in view of highlighting the non-invasive feature of our strategy.

The paper is organized as follows: Section~\ref{sec:link_IGA_FEM} reviews the existing FEM-to-IGA link which constitutes the necessary prerequisite for a seamless coupling between global-IGA and local-FEM. Then, Section~\ref{sec:non_inv_iga_fem} is devoted to the derivation of our fully non-invasive global-IGA/local-FEM methodology. Finally, in Section~\ref{sec:results}, the performance of the proposed implementation is demonstrated through a series of benchmarks involving complex local behaviors, before concluding remarks are formulated in Section~\ref{sec:conclu}.

\section{IGA as a projection of FEM onto a regular reduced basis}
\label{sec:link_IGA_FEM}

This section undertakes to outline IGA as a projection of FEM onto a regular reduced basis. An alternative lighting on the relation between IGA and FEM is thus provided compared to the more common one that consists of viewing IGA as encompassing FEM. This leads to formulate a complete algebraic bridge that directly goes from Lagrange nodal polynomials to B-Spline and NURBS functions, which builds the foundations to achieve our fully non-invasive global-IGA/local-FEM coupling.

\subsection{Some necessary ingredients for B-Splines and NURBS}
\label{sub:spline_base}

Let us start by recalling key elements regarding B-Spline and NURBS~\cite{Cohen80,Piegl97,Massarwi19}. Only the fundamentals are outlined in the following. For further details, besides the pioneering contributions~\cite{hughes05,cottrell09}, the reader is referred to the works cited hereafter.

The NURBS functions lend themselves to an exact representation of many shapes used in engineering, such as conical sections (circles, cylinders, etc). A general expression for a NURBS geometry $G^h(\xi)$ with parameter $\xi \in \mathbb{R}^d$ can be written as:
\begin{equation}
G^h(\xi) = \sum_{i=1}^{n_{IG}} R_i (\xi) x_i^{IG} = \mathbf{R}^T(\xi)\mathbf{x}^{IG},
\label{eq:geo_nurbs}
\end{equation}
where $\mathbf{R}$ and $\mathbf{x}^{IG}$ denote the matrix of the $n_{IG}$ NURBS basis functions and the vector collecting the locations of the associated control points, respectively. A NURBS entity in $\mathbb{R}^d$ can be seen as a projection of a B-Spline entity in $\mathbb{R}^{d+1}$, which results in expressing the NURBS functions from the B-Spline ones by associating a weight to each control point (see, \emph{e.g.},~\cite{hughes05,cottrell09}). Then, all one needs to do in order to define the multivariate B-Spline function $N_i$ at control point $i$ is to perform the tensor product of the univariate B-Spline functions associated with this point in the different parametric directions. In the end, the $n^1_{IG}$ univariate B-Spline basis functions are piecewise polynomials defined by their polynomial degree $p$ and a set of non-decreasing parametric coordinates $\xi_1^i \in \mathbb{R}$ collected into a knot-vector $\Xi_1 = \left\lbrace\xi_1^1,\xi_1^2,..,\xi_1^{n^1_{IG}+p+1} \right\rbrace$.

\begin{figure}[htb]{
		\centering
		\includegraphics[width=1\linewidth]{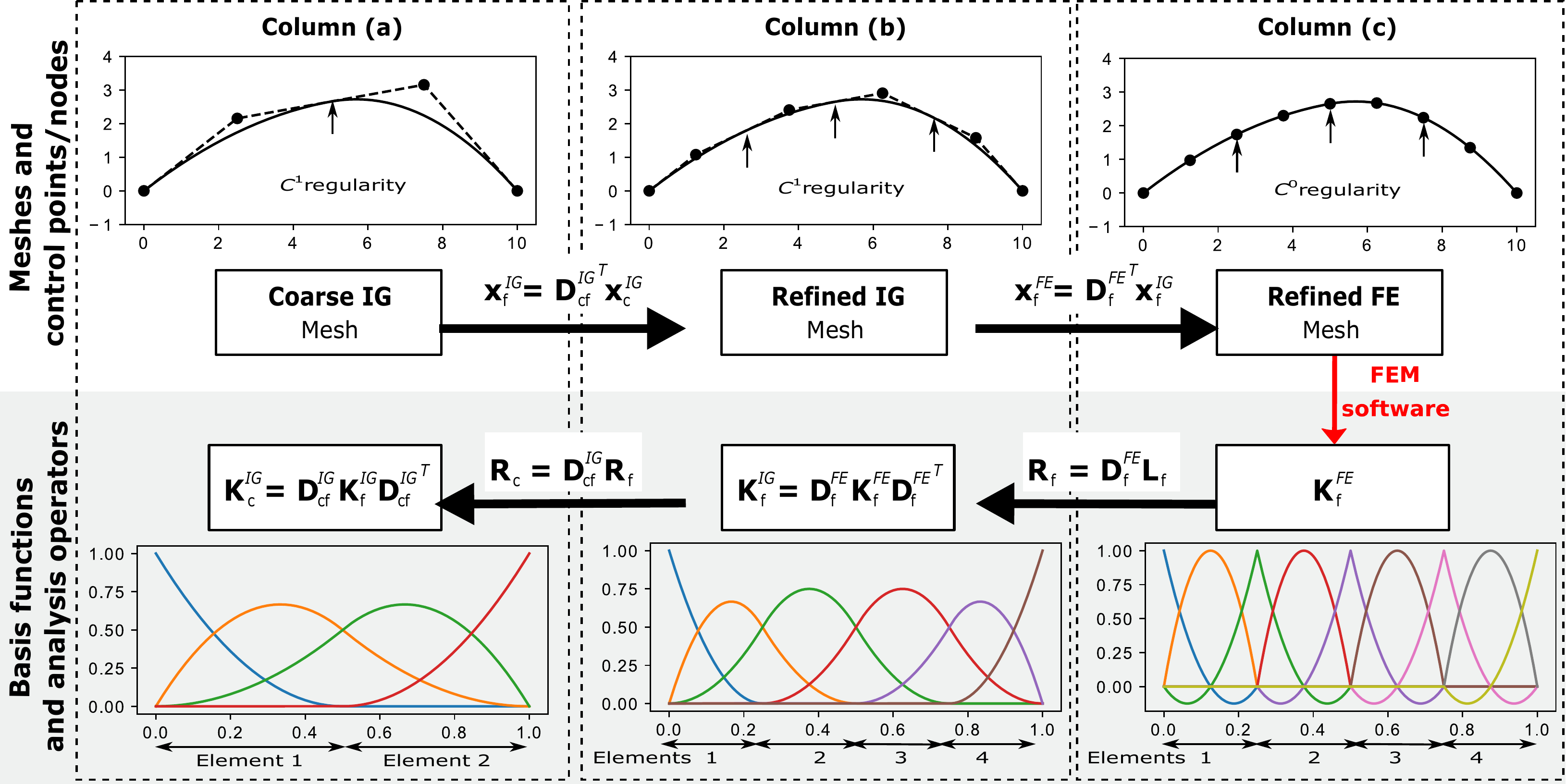}
		\caption{Link between IGA and FEM. Column (a):  the initial IG discretization (four control points associated with four global quadratic spline functions), column (b): the refined IG discretization (knots 0.25 and 0.75 are added which results in the definition of 6 control points to maintain the geometry), column (c):  the refined FE discretization (knots 0.25, 0.5 and 0.75 are added to obtain a $C^0$ regularity at those locations and a Lagrange-to-Bernstein change of basis is performed which results in the construction of 9 FE nodes to recover the geometry). The different linear operators can then be used to recover, by projection, the refined and initial IG stiffness matrices from the refined FE one computed using a classical FE software, taken as a black- box.}
		\label{fig:total_iga_fem_link}}
\end{figure}

The interesting feature of splines is their higher degree of regularity. Indeed, a B-Spline function of degree $p$ can reach a $C^{p-1}$ regularity at knot $\xi_1^i$ if this one is single in $\Xi_1$.  As an illustration, let us refer to Fig.~\ref{fig:total_iga_fem_link} (in particular, see bottom in columns (b) and (c)). The six global B-Spline functions attain a $C^1$-regularity at knots 0.25, 0.5 and 0.75 while the FE space is built from nine global shape functions that meet a $C^0$-regularity at those locations. Consequently, for a given polynomial degree and a similar number of elements, a $C^{p-1}$ B-Spline mesh comes with less DOF than the corresponding $C^0$ FE mesh, which is totally understandable since the space of $C^{p-1}$ functions is included into the space of $C^0$ functions. For instance, for a 3D solid mesh composed of $100$ elements in the three parametric directions: we gain a factor of about $8$, $25$ and $57$ for $p=2$, 3 and 4, respectively. This is the main feature of IGA that provides increased per-DOF accuracy with respect to FEM when smooth solutions are to be captured. IGA will thus be used in this work for the global model, as regular solutions are expected at this level.

Finally, spline functions present efficient smooth refinement procedures which allow to enhance the approximation space without changing the geometry.  An example is given in Fig.~\ref{fig:total_iga_fem_link} (columns (a) and (b)). In practice, matrix representation of the spline refinement procedures exist~\cite{Piegl97}; that is, denoting by $\mathbf{R}_{\mathrm{c}}$ and $\mathbf{R}_{\mathrm{f}}$ (resp.  $\mathbf{x}^{IG}_{\mathrm{c}}$ and $\mathbf{x}^{IG}_{\mathrm{f}}$), the matrices (resp. vectors) collecting the coarse and fine spline functions (resp. control points), we can build the refinement operator $\mathbf{D}^{IG}_{\mathrm{cf}}$ such that:
\begin{equation}
\mathbf{R}_{\mathrm{c}}= \mathbf{D}^{IG}_{\mathrm{cf}}~\mathbf{R}_{\mathrm{f}} \qquad \mathrm{and} \qquad \mathbf{x}^{IG}_{\mathrm{f}}= \left(\mathbf{D}^{IG}_{\mathrm{cf}}\right)^T\mathbf{x}^{IG}_{\mathrm{c}}.
\label{eq:NcCrNf}
\end{equation}
For more details on refinement strategies of splines and their matrix representations, reference is made to~\cite{Piegl97,Cottrell07,Bouclier_book_22}.

\subsection{The link between IGA and FEM}
\label{sub:iga_fem_link}

The objective now is to relate IGA and FEM.  In order to do so, we make use of previous works~\cite{Tirvaudey19,Colantonio20}  in which a global algebraic bridge between IGA and FEM was established by resorting to Bézier-based operators~\cite{Borden11,Schillinger16,Kamensky19}. Although this is not requested for our hybrid global-IGA/local-FEM scheme, we also show how to build, in an explicit and plug-and-play manner, the IG operators (stiffness matrix and load vector) from their FE counterparts computed using a standard FE code. This highlights our point of view on IGA (projection of FEM onto a more regular, reduced basis) and will allow to perform the implementation even without having a global IG code in hand, as will be shown in Section~\ref{sec:non_inv_iga_fem} (see remark~\ref{rk:FEM/FEM}). Finally, we refer again to Fig.~\ref{fig:total_iga_fem_link} to illustrate the theory on a concrete and simple example.

\subsubsection{Viewpoint of the technologies}

Starting with B-Spline functions, it is easy to build $C^0$ polynomials: it suffices to repeat all the interior knots of the knot-vectors until they reach a $p$ multiplicity. This treatment actually consists in the Bézier extraction~\cite{Borden11}, which leads to the construction of Bernstein polynomials. The advantage of Bernstein functions is that they exhibit an elementary structure which is similar to FEM. Thus, to arrive at nodal Lagrange polynomials, the remaining task is to perform a change of basis to scale the basis functions so that they are equal to one at the corresponding FE nodes. The composition of the Bézier extraction with the Lagrange-to-Bernstein change of basis yields the Lagrange extraction~\cite{Schillinger16,Tirvaudey19} that allows  to formulate a smooth polynomial B-Spline discretization in terms of a standard FE discretizations. The Lagrange extraction is illustrated in Fig.~\ref{fig:total_iga_fem_link} for a quadratic spline curve made of four elements (see columns (b) and (c)). The knots 0.25, 0.5 and 0.75 are added to first apply the Bézier extraction and a Lagrange-to-Bernstein change of basis is then performed to reach standard nodal Lagrange polynomials.

The Lagrange extraction can be extended to the case of NURBS: it consists in expressing NURBS functions in terms of rational Lagrange functions, as detailed in~\cite{Schillinger16}. 
To truly involve Lagrange polynomials, a projection allowing to formulate rational functions in terms of polynomials is required. In order to do so, we proceed as in~\cite{Tirvaudey19}; that is, we start with the NURBS version of the Lagrange extraction and then perform the projection at the Lagrange level by approximating to $1$ all the weights of the rational Lagrange functions. In the end, we are able to build refined B-Spline or NURBS basis functions $\mathbf{R}_{\mathrm{f}}$ from standard refined Lagrange polynomials $\mathbf{L}_{\mathrm{f}}$:
\begin{equation}
\mathbf{R}_{\mathrm{f}}= \mathbf{D}^{FE}_{\mathrm{f}}~\mathbf{L}_{\mathrm{f}},
\label{eq:NftoLf}
\end{equation}
where $\mathbf{D}^{FE}_{\mathrm{f}}$ is a linear operator that traduces the FEM-to-IGA bridge. In addition, we can build a FE mesh that represents the same geometry as the spline one by taking the FE nodes $\mathbf{x}^{FE}_{\mathrm{f}}$ such that:
\begin{equation}
\mathbf{x}^{FE}_{\mathrm{f}}= \left(\mathbf{D}^{FE}_{\mathrm{f}}\right)^T\mathbf{x}^{IG}_{\mathrm{f}}.
\label{eq:fe_nodes}
\end{equation}
For more details on this topic, the interested reader can refer to~\cite{Tirvaudey19}.

\begin{remark}
	Let us underline that equalities~\eqref{eq:NftoLf} and~\eqref{eq:fe_nodes} do not strictly hold in case of NURBS. Indeed, moving from a rational to a polynomial geometry necessarily leads to some approximations. However, as demonstrated in~\cite{Tirvaudey19}, the error related to such approximations is largely insignificant when considering refined geometries compared to the associated NURBS discretization error. The results of the present paper in section~\ref{sec:results} will also confirm this statement.
\end{remark}

\subsubsection{Viewpoint of the resulting approximation spaces}

From an analysis point of view, matrix $\left(\mathbf{D}^{FE}_{\mathrm{f}}\right)^T$ can be seen as a collection of modes, each mode being a column of the matrix. Considering the displacement as the field of interest, each of these modes is the displacement of all FE nodes caused by a unitary displacement of a control point. This collection of modes is a basis for a vector subspace of the vector space generated by the FE functions. As a result, if one does not have an IG code in hand, one can construct a FE mesh from the IG mesh in a pre-processing step (see Fig.~\ref{fig:total_iga_fem_link}(top)), and then simply perform model reduction~\cite{Quarteroni15,Hesthaven16,Benner22} to obtain the IG linear system to solve from the FE one.  In accordance with the Ritz-Galerkin method, the reduced-order problem writes:
\begin{equation}
\mathbf{K}_{\mathrm{f}}^{IG}\mathbf{u}_{\mathrm{f}}^{IG} = \mathbf{f}_{\mathrm{f}}^{IG} \quad \Leftrightarrow \quad \mathbf{D}^{FE}_{\mathrm{f}}\mathbf{K}_{\mathrm{f}}^{FE}\left(\mathbf{D}^{FE}_{\mathrm{f}}\right)^T \mathbf{u}_{\mathrm{f}}^{IG} = \mathbf{D}^{FE}_{\mathrm{f}} \mathbf{f}_{\mathrm{f}}^{FE},
\end{equation}
where $\mathbf{K}_{\mathrm{f}}^{IG}$ (respectively $\mathbf{f}_{\mathrm{f}}^{IG}$) is the stiffness matrix (resp. load vector) associated with the refined IG mesh (see Fig.~\ref{fig:total_iga_fem_link}(column(b)), and $\mathbf{K}_{\mathrm{f}}^{FE}$ (resp. $\mathbf{f}_{\mathrm{f}}^{FE }$) is the stiffness matrix (resp. load vector) associated with the refined FE mesh (see Fig.~\ref{fig:total_iga_fem_link}(column(c)). This enables to compute the refined IG DOF vector $\mathbf{u}_{\mathrm{f}}^{IG}$ in a non-invasive manner from FEM: it does not require any modifications in how $\mathbf{K}_{\mathrm{f}}^{FE}$ and $\mathbf{f}_{\mathrm{f}}^{FE }$ are built from a standard FE code. Then, note that the resulting IG displacement can be back-converted in terms of nodal displacements:
\begin{equation}
\mathbf{u}^{FE}_{\mathrm{f}}= \left(\mathbf{D}^{FE}_{\mathrm{f}}\right)^T\mathbf{u}^{IG}_{\mathrm{f}},
\label{eq:link_disp}
\end{equation}
so that existing subroutines of the FE code can be used for post-processing. Similarly, the initial coarse IG solution $\mathbf{u}_{\mathrm{c}}^{IG} $ can be obtained by composing with the spline refinement operator $\mathbf{D}^{IG}_{\mathrm{cf}}$:
\begin{equation}
\begin{split}
\mathbf{K}_{\mathrm{c}}^{IG}\mathbf{u}_{\mathrm{c}}^{IG} = \mathbf{f}_{\mathrm{c}}^{IG} \quad \Leftrightarrow \quad &\mathbf{D}^{IGFE}_{\mathrm{cf}}\mathbf{K}_{\mathrm{f}}^{FE}\left(\mathbf{D}^{IGFE}_{\mathrm{cf}}\right)^T \mathbf{u}_{\mathrm{c}}^{IG} = \mathbf{D}^{IGFE}_{\mathrm{cf}} \mathbf{f}_{\mathrm{f}}^{IG}, \\ \\
&\mathrm{with} \quad \mathbf{D}^{IGFE}_{\mathrm{cf}}=\mathbf{D}^{IG}_{\mathrm{cf}}\mathbf{D}^{FE}_{\mathrm{f}}.
\end{split}
\label{eq:reduced_order}
\end{equation}
Once again we refer to Fig.~\ref{fig:total_iga_fem_link} that gives an overview of the different transformations allowing to perform IGA from FEM and to~\cite{Tirvaudey19,Colantonio20} for further details on this topic.

\begin{remark}
	From this point of view, IGA consists of nothing more than projecting FEM onto a reduced, more regular vector subspace made of spline functions. Note however that in contrast to more standard reduced basis methods, operator $\mathbf{D}^{FE}_{\mathrm{f}}$ can still appear quite large but it is highly sparse here.
\end{remark}

\begin{remark}
	Let us also underline that the concept of Bézier extraction has now been generalized to a large variety of advanced splines (such as T-Splines~\cite{Scott11,Kamensky19},  hierarchical B-Splines and NURBS~\cite{Hennig16,Angella20}, and hierarchical T-Splines~\cite{Evans15,Chen18}). As a result, the developed FEM-to-IGA bridge and its use for achieving a non-invasive implementation with respect to FEM could be straightforwardly extended to many other spline technologies.
\end{remark}


\section{Non-invasive global-IGA/local-FEM}
\label{sec:non_inv_iga_fem}

Now that the link between IGA and FEM has been reviewed, let us introduce our hybrid global-IGA/local-FEM scheme that seeks to combine the interests of both technologies for global/local simulations. Let us underline at this stage that our approach is generic in terms of programming environments: the users may have in hand an IG code (performing standard elasticity) and wish to couple it with a specific FE software to model complex local phenomena, or the users only have at their disposal FE packages. The aim is to arrive at an automatic coupling between global-IGA and local-FEM.  For the presentation, the non-invasive global/local algorithm is first briefly given and then, the specific construction of the models along with the dedicated implementation to achieve a fully non-invasive hybrid IGA/FEM procedure is detailed.

\subsection{The reference algorithm}
\label{sub:non_inv_ref}

\subsubsection{The global/local problem}

\begin{figure}[htb]{
		\centering
		\includegraphics[width=\linewidth]{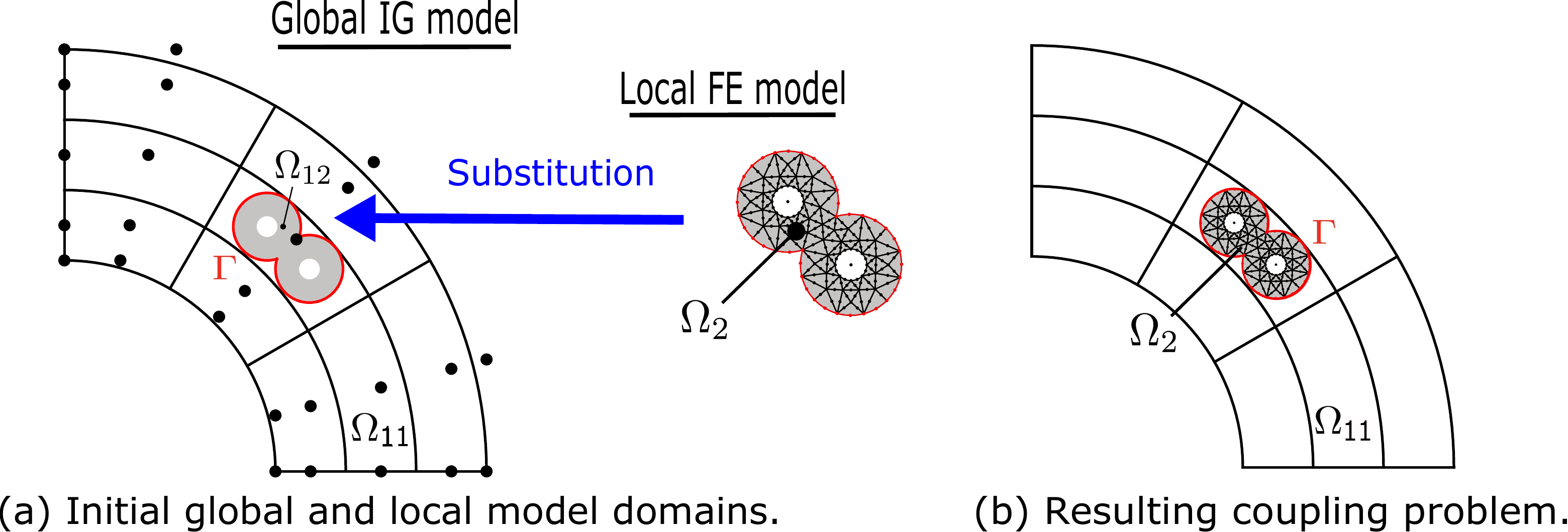}
		\caption{Example of a global-IGA/local-FEM problem. The global IG model over sub-domain $\Omega_{12}$ is replaced by the finer local FE model of domain $\Omega_{2}$ through interface $\Gamma$, which enables to integrate geometrical details (holes) along with possible non-linearities around (\emph{e.g.}, contact, cracks, plasticity) within the initial regular coarse model.}
		\label{fig:global_local_non_inv}}
\end{figure}

Let us first properly introduce the mechanical multiscale global/local problem that we seek to solve. We consider a global (coarse) IG model of a structure that is characterized by a physical domain $\Omega_{1}$ (see Fig.~\ref{fig:global_local_non_inv}(a)(left)). This domain is divided into two disjoint, open and bounded subsets $\Omega_{11}$ and $\Omega_{12}$. Those two non-overlapping sub-domains share a common interface denoted by $\Gamma$ such that $\Omega_1 = \Omega_{11} \cup \Omega_{12} \cup \Gamma$ and $\Omega_{11} \cap \Omega_{12}=\varnothing$. A simple linear elastic modeling is adopted for the global structure. We assume that such a behavior and the coarse spline discretization is sufficient to accurately capture the solution except in the small region $\Omega_{12}$ where a local (possibly non-smooth, singular, or even discontinuous) phenomenon is to be introduced. As a consequence, a local, more detailed FE "sub-model" characterized by domain $ \Omega_{2}$ is constructed to replace the global model in $\Omega_{12}$ (see Fig.~\ref{fig:global_local_non_inv}(a)(right)). The substitution of the FE local model within the IG global one is achieved through interface $\Gamma$. The resulting global/local problem to be solved is a hybrid IGA/FEM multi-domain problem in $\Omega_{11} \cup \Omega_{2} \cup \Gamma$, the global solution in $\Omega_{12}$ being discarded (see Fig.~\ref{fig:global_local_non_inv}(b)).

Although the standard non-invasive method applies for any (possibly non-linear) local behavior, we will consider for simplicity in the presentation in this section that the local model is also linear elastic. Yet non-linear local models will be investigated in the numerical results section~\ref{sec:results}.

\subsubsection{Monolithic solution}

The usual starting point in the derivation of the non-invasive global/local coupling strategy is to weakly formulate the coupling with a Lagrange multiplier approach (also called Mortar approach, see~\cite{Brivadis15,Bouclier17,Miao20} among others). We directly write below the formulation in the discrete setting. We denote the matrices of shape functions associated with the global IG model over $\Omega_1$ and local FE model over $\Omega_2$ by $\mathbf{R}_{\mathrm{{1c}}}$ and $\mathbf{L}_2$, respectively. We also need to define $\mathbf{R}_{\mathrm{{11c}}}$ that stands for the matrix that gathers the restricted part to sub-domain $\Omega_{11}$ of the shape functions of the global model. The displacement DOF vectors corresponding to each of the models are $\mathbf{u}_{\mathrm{1c}}^{IG}$, $\mathbf{u}_{\mathrm{2}}^{FE}$ and $\mathbf{u}_{\mathrm{11c}}^{IG}$, respectively. Then, a Lagrange mutliplier field defined on $\Gamma$ is introduced, as a dual unknown, to represent the interface traction forces on both sides of the interface. The associated DOF vector is denoted by $\boldsymbol{\lambda}$. In addition,  let us write $\mathbf{\Phi}_{\Gamma}$ at this stage for the shape function matrix corresponding to the Lagrange multiplier field. The  construction of $\mathbf{\Phi}_{\Gamma}$ will be given further (see section~\ref{sub:mesh_build_non_inv}).
With all the above notations, the Lagrange multiplier approach for the resulting coupling problem (see Fig.~\ref{fig:global_local_non_inv}(b)) finally leads to solving:
\begin{equation}
\def\arraystretch{1.5}
\begin{bmatrix}
\mathbf{K}_\mathrm{11c}^{IG} & \mathbf{0}  & {\mathbf{C}_\mathrm{11c}^{IG}}^T  \\
\mathbf{0} & \mathbf{K}_2^{FE} & -{\mathbf{C}_2^{FE}}^T  \\
\mathbf{C}_\mathrm{11c}^{IG} & -\mathbf{C}_2^{FE} & \mathbf{0} \\
\end{bmatrix}
\begin{pmatrix}
\mathbf{u}_\mathrm{11c}^{IG} \\  \mathbf{u}_2^{FE}  \\ \boldsymbol{\lambda} 
\end{pmatrix}
=
\begin{pmatrix}
\mathbf{f}_\mathrm{11c}^{IG}  \\  \mathbf{f}_2^{FE} \ \\  \mathbf{0} 
\end{pmatrix},
\label{eq:syst_disc_ref}
\end{equation}
where $\mathbf{K}_\mathrm{11c}^{IG}$ (respectively $\mathbf{f}_\mathrm{11c}^{IG}$) and $\mathbf{K}_2^{FE}$ (resp. $\mathbf{f}_2^{FE}$) are the classical stiffness matrices (resp. load vectors) associated with sub-domains $\Omega_{11}$ and $\Omega_2$, and $\mathbf{C}_\mathrm{11c}^{IG}$ and $\mathbf{C}_2^{FE}$ are the Mortar coupling operators that formally read:
\begin{equation}
\mathbf{C}_\mathrm{11c}^{IG} = \int_{\Gamma} \mathbf{\Phi}_\Gamma \mathbf{R}_{\mathrm{{11c}}}^T\mathrm{d}\Gamma \quad ; \quad
\mathbf{C}_2^{FE} = \int_{\Gamma} \mathbf{\Phi}_\Gamma \mathbf{L}_2^T\mathrm{d}\Gamma.
\label{eq:expr_mortar}
\end{equation}

The resolution~\eqref{eq:syst_disc_ref} of the global/local problem constitutes the classical monolithic approach: the  coupled model of Fig.~\ref{fig:global_local_non_inv}(b) is computed directly using a single direct solver. This strategy is invasive in the sense that it requires (i) to modify the initial global model (and thus its operators) to remove some of its elements, or possibly pieces of elements (see Fig.~\ref{fig:global_local_non_inv} again) which may lead to ill-conditioned stiffness operators~\cite{Prenter17,Wei21}, and (ii) to set up an additional solver that merges the contributions of the two models. In case the local detail grows up (during crack propagation, or expansion of damage or plasticity for instance), the situation is getting even worse since not only the local operator $\mathbf{K}_{2}$ but also the global operator $\mathbf{K}_{11}$ have to be fully re-built, and the augmented system~\eqref{eq:syst_disc_ref} re-factorized during the simulation.

\subsubsection{Non-invasive iterative solution}

\begin{figure}[htb]{
		\centering
		\includegraphics[width=\linewidth]{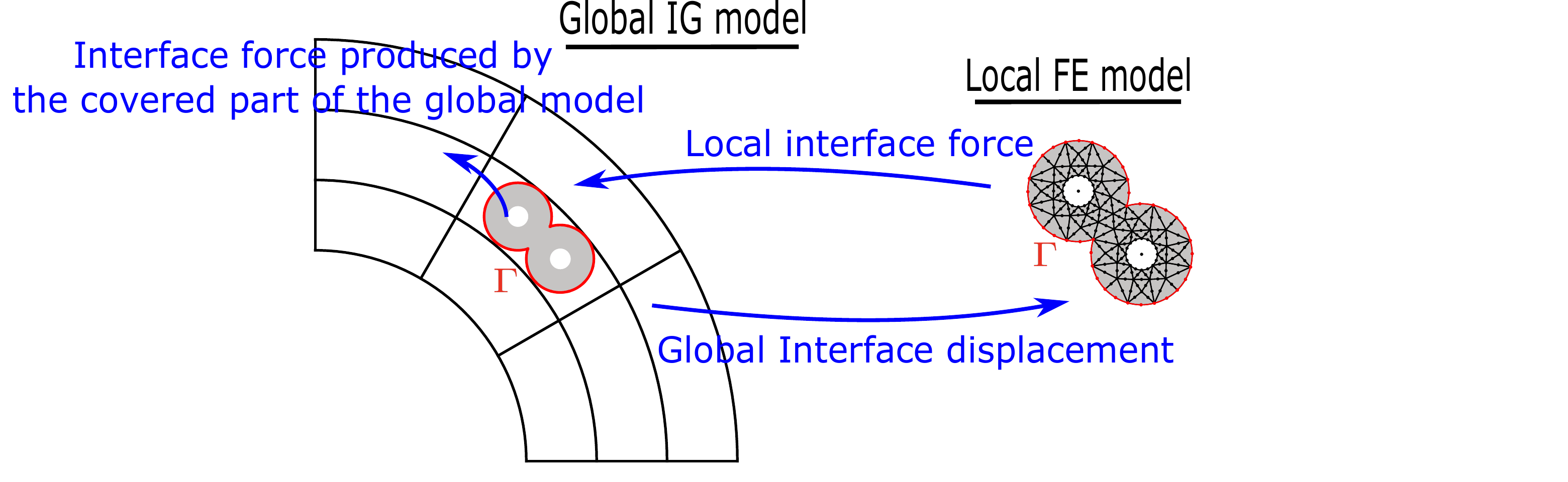}
		\caption{The iterative non-invasive exchange procedure.  Starting with a global Neumann resolution (with respect to $\Gamma$), the interface displacement is transferred from the global to the local model. Then, a local Dirichlet problem is solved and the interface traction force coming from the local model is applied to the global model along with the interface traction force produced by the covered part of the global model at previous iteration.}
		\label{fig:non_inv_classic}}
\end{figure}

Conversely, the non-invasive strategy is based on an iterative exchange procedure that alternates between global solutions over $\Omega_1$ and local solutions over $\Omega_{2}$. Briefly, its derivation is performed in two steps. First, we split the initial system~\eqref{eq:syst_disc_ref} in order to identify (in terms of boundary conditions applied on $\Gamma$) a Neumann and a Dirichlet problem over $\Omega_{11}$ and $\Omega_{2}$, respectively. Then,  we make use of the available continuous prolongation of the global solution from $\Omega_{11}$ to $\Omega_{12}$. This allows to apply the additivity of the integral with respect to domain $\Omega_1 = \Omega_{11} \cup \Omega_{12} \cup \Gamma$ to recover the initial whole global model. In the end, we obtain the following asymmetric algorithm in the sense that Dirichlet and Neumann problems with respect to $\Gamma$ are alternatively solved until convergence. More precisely, for the $n$th iteration, starting with initial guesses $\boldsymbol{\lambda}^{(0)}$ and ${\mathbf{u}_\mathrm{1c}^{IG}}^{(0)}$, we look for ${\mathbf{u}_\mathrm{1c}^{IG}}^{(n)}$, ${\mathbf{u}_2^{FE}}^{(n)}$ and $\boldsymbol{\lambda}^{(n)}$ such that:
\begin{enumerate}
	\item Resolution of a Neumann problem (with respect to $\Gamma$) over $\Omega_{1}$:
	\begin{equation}
	\def\arraystretch{1.5}
	\mathbf{K}_\mathrm{1c}^{IG} ~{\mathbf{u}_\mathrm{1c}^{IG}}^{(n)} = \mathbf{f}_\mathrm{1c}^{IG}  -{\mathbf{C}_\mathrm{1c}^{IG}}^T \boldsymbol{\lambda}^{(n-1)}+{\textcolor{red}{\overline{\boldsymbol{\lambda}}}_\mathrm{12c}^{IG}}^{(n-1)}.
	\label{eq:nonint_global_dis}
	\end{equation}
	\item Resolution of a Dirichlet problem (with respect to $\Gamma$) over $\Omega_{2}$:
	\begin{equation}
	\def\arraystretch{1.5}
	\begin{bmatrix}
	\mathbf{K}_{2}^{FE} & -{\mathbf{C}_2^{FE}}^T  \\
	-\mathbf{C}_2^{FE} & \mathbf{0}  \\
	\end{bmatrix}
	\begin{pmatrix}
	{\mathbf{u}_2^{FE}}^{(n)} \\ \boldsymbol{\lambda}^{(n)} 
	\end{pmatrix}
	=
	\begin{pmatrix}
	\mathbf{f}_2^{FE} \\  -\mathbf{C}_\mathrm{1c}^{IG} {\mathbf{u}_{1}^{IG}}^{(n)}   \end{pmatrix}.
	\label{eq:nonint_local_dis}
	\end{equation}
\end{enumerate}
In the above equations, $\mathbf{K}_\mathrm{1c}^{IG}$, $\mathbf{f}_\mathrm{1c}^{IG}$ and $\mathbf{C}_\mathrm{1c}^{IG}$ simply consist of the prolongation of former operators $\mathbf{K}_\mathrm{11c}^{IG}$, $\mathbf{f}_\mathrm{11c}^{IG}$ and $\mathbf{C}_\mathrm{11c}^{IG}$, respectively, from $\Omega_{11}$ to $\Omega_1$. $\textcolor{red}{\overline{\boldsymbol{\lambda}}}_\mathrm{12c}^{IG}$ is introduced to denote the discrete reaction forces at $\Gamma$ produced by the covered part $\Omega_{12}$ of the global model. It emerges to counterbalance the effect of this covered region since this one is not present in the reference coupling problem (see Fig.~\ref{fig:global_local_non_inv}(b) again). It reads at iteration $n-1$:
\begin{equation}
{\overline{\boldsymbol{\lambda}}_\mathrm{12c}^{IG}}^{(n-1)} = \overline{\mathbf{K}}_\mathrm{12c}^{IG} {\mathbf{u}_\mathrm{1c}^{IG}}^{(n-1)}- \overline{\mathbf{f}}_\mathrm{12c}^{IG},
\label{eq:R12}
\end{equation}
where $\overline{\mathbf{K}}_\mathrm{12c}^{IG}$ and $\overline{\mathbf{f}}_\mathrm{12c}^{IG}$ are the extensions to $\Omega_1$ of the classical stiffness matrix $\mathbf{K}_\mathrm{12c}^{IG}$ and load vector $\mathbf{f}_{12c}^{IG}$ of $\Omega_{12}$, respectively. They literally contain the classical stiffness and load vector operators for $\Omega_{12}$ and are padded with zeros to make them the same dimension of $\mathbf{u}_\mathrm{1c}^{{IG}^{(n)}}$. The procedure can be interpreted as a fixed point strategy aiming at ensuring the equilibrium of the interface reaction forces (see Eq.~\eqref{eq:nonint_global_dis}) provided that the displacement is transferred at each iteration between the two models (see Eq.~\eqref{eq:nonint_local_dis}). An illustration of the algorithm is provided in Fig.~\ref{fig:non_inv_classic}. For more information, we advise the interested reader to consult the following reviews on the subject~\cite{Duval16,Allix20}.

Algorithm~\eqref{eq:nonint_global_dis}-\eqref{eq:nonint_local_dis} constitutes what is commonly referenced as the non-invasive global/local coupling strategy in the literature. Since the initial global IG model is now unmodified, its stiffness operator remains well-conditioned regardless of the shape of the local region, and it can be assembled and factorized only once during the pre-processing step~\cite{Bouclier16,Bouclier18}. Furthermore, the global and local problems being solved alternatively and the interaction between the two models being restricted to interface $\Gamma$, the formalism offers the possibility to couple an IG and a FE code with very few implementation effort. The price to pay is the number of iterations of the fixed point solver. However, this one can be deeply reduced by means of accelerations techniques, such as based on an Aitken's Delta Squared method or a Quasi-Newton method (see~\cite{Duval16,Gosselet18} to name a few). Even more important, such acceleration techniques may appear necessary to ensure the convergence of the algorithm in challenging situations (see, \emph{e.g.},~\cite{Chevreuil13} regarding the theory).

\begin{remark}
	We recall that the method is by no means limited to a linear elastic local model. Indeed, as long as we are able to apply Dirichlet boundary conditions to the local problem and to compute (directly or in a post-processing step) the corresponding reactions forces, any local behavior can be considered, as demonstrated in the large literature on the topic (see~\cite{Gendre09,Duval16,Guinard18,Duarte21}  to name a few), and as will be shown later in this paper in Section~\ref{sec:results}.
\end{remark}

\subsection{Construction of the FE model to reach a conforming global/local interface}
\label{sub:mesh_build_non_inv}

Incorporating a specific local region in an IG patch without care may result in the overlap of some global knot-span elements due to the rigid tensor product structure of (standard) multivariate spline bases (see Figs.~\ref{fig:global_local_non_inv}(a) and~\ref{fig:non_inv_classic} again). More precisely, the difficulty relies on (i) the evaluation of integrals over pieces of knot-span elements (to get the interface reaction force $\textcolor{red}{\overline{\boldsymbol{\lambda}}}_\mathrm{12c}^{IG}$, see Eq.~\eqref{eq:R12}), and (ii) the formulation of a coupling method adapted to an immersed interface. As a remedy, the idea here is to call upon efficient (classic) FE meshing procedures to reach a conforming interface (similar pragmatic approach as in~\cite{Guinard18} in the context of global/local FEM). In order to do so in a simple, automatic and consistent way, we make use of the FEM-to-IGA bridge of Section~\ref{sec:link_IGA_FEM}, which will also enable to arrive at a fully non-invasive strategy in the sense that not only the global/local coupling is non-invasive but also the construction of the coupling operators (see Eq.~\eqref{eq:expr_mortar}) from only FE resources.

\begin{figure}[htb]{
		\centering
		\includegraphics[width=\linewidth]{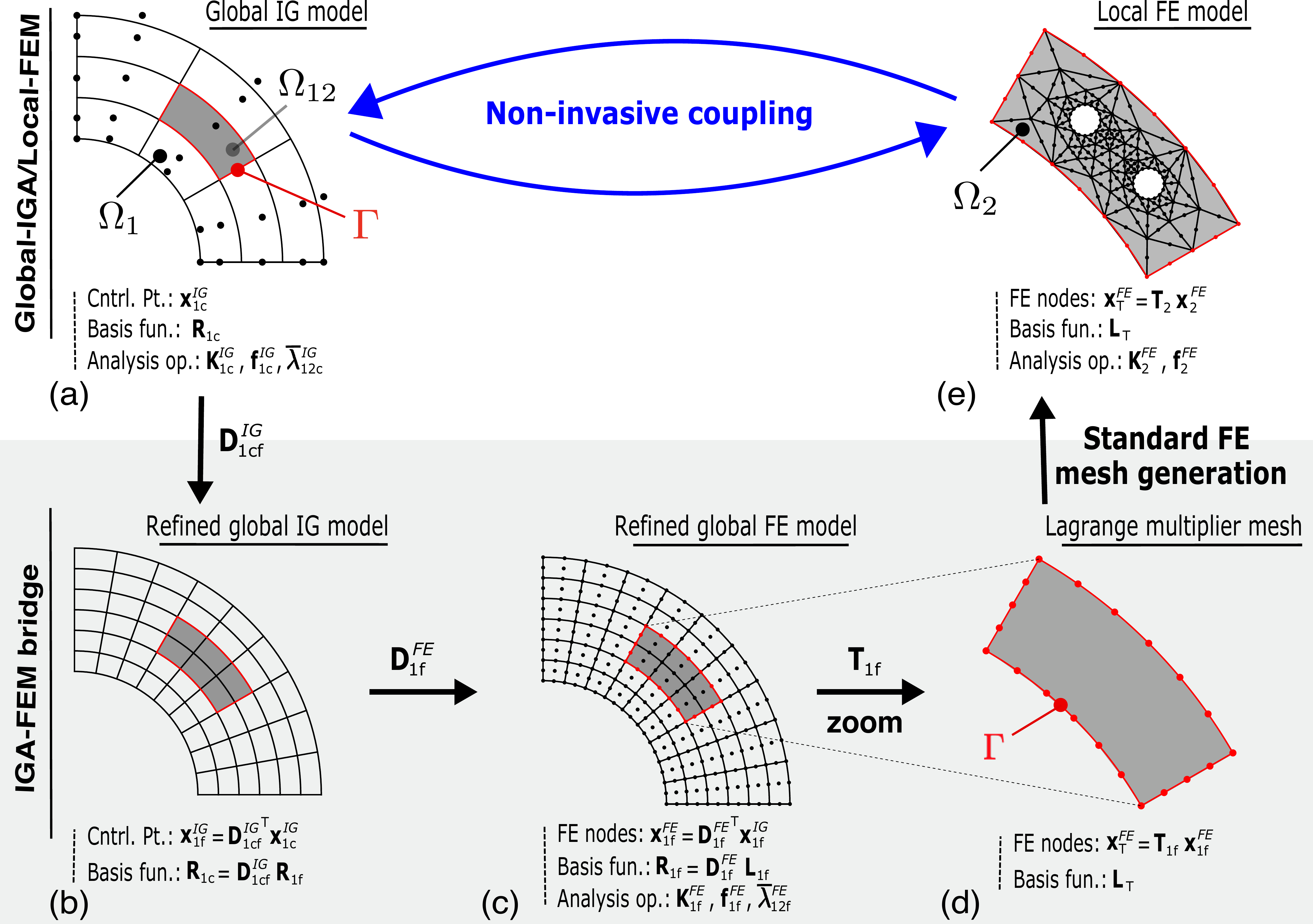}
		\caption{Illustration of the proposed procedure to build in a simple, automatic and consistent way a conforming global-IGA/local-FEM discretization. The strategy strongly relies on the FEM-to-IGA bridge, which also offers the opportunity to simplify the implementation of the mechanical solver by involving only  FE operators.}
		\label{fig:gIGAlFEM_papier_principe}}
\end{figure}

The proposed procedure for the construction of the conforming global-IGA/local-FEM modeling is illustrated in Fig.~\ref{fig:gIGAlFEM_papier_principe}. This figure also presents the notations followed which are consistent with  all those introduced previously. Starting with a global  IG model of the whole structure (a), a specific local FE model (e), meant to replace the global IG model in an area of interest, is built. The area of interest is obviously chosen as a set of initial global knot-span elements so that we are able to reach a conforming coupling interface. More precisely, from the initial global IG model (a), we apply standard spline refinement procedures~\eqref{eq:NcCrNf} (associated operator $\mathbf{D}^{IG}_{\mathrm{1cf}}$) to obtain the refined global IG model (b). Then, we make use of the FEM-to-IGA bridge~\eqref{eq:fe_nodes} (operator $\mathbf{D}^{FE}_{\mathrm{1f}}$) to obtain the corresponding refined global FE model (c). It is therefore possible to extract the FE description of the interface (d) by calling upon a trace operator ($\mathbf{T}_{\mathrm{1f}}$) that selects only the nodes (or the DOF) concerned with the interface. Finally, by applying existing optimized FE meshing procedures (\emph{e.g.}, GMSH~\cite{GMSH}, or Salome-Meca~\cite{salome} which will be used for the numerical illustrations), we can build a local FE mesh that is conforming with the interface. In order to do so, it may be noticed at this stage that the constructed FE mesh must have the same polynomial degree as the interface. 

As is standard practice in conforming global/local FEM, we choose to take the trace along the interface of the functions of the local model ($\mathbf{L}_{\mathrm{T}}$) to discretize the Lagrange multiplier field. We thus now replace $\mathbf{\Phi}_{\Gamma}$ by $\mathbf{L}_{\mathrm{T}}$  in~\eqref{eq:expr_mortar} which leads to :
\begin{equation}
\mathbf{C}_2^{FE} = \int_{\Gamma} \mathbf{L}_{\mathrm{T}} \mathbf{L}_{\mathrm{2}}^T \mathrm{d} \Gamma,
\end{equation}
\emph{i.e.}, the mass matrix associated with the local FE model of the FE interface. We never encountered instabilites in our numerical experiments with such a choice.

\subsection{Fully non-invasive implementation of global-IGA/local-FEM}
\label{sub:fully_non_inv}

In general, the computation of the Mortar coupling operators may not appear trivial. It is necessary to build an integration technique on the interface and to have access to the values of the shape functions on both sides at each integration point, which is not a classical output of industrial codes and breaks the concept of non-invasiveness. Here the proposed meshing strategy depicted in Fig.~\ref{fig:gIGAlFEM_papier_principe} allows to circumvent the difficulty; thanks to the FEM-to-IGA bridge, the hybrid IGA/FEM coupling operators become explicit. Indeed,  we can compute $\mathbf{C}_\mathrm{1c}^{IG}$ as follows:
\begin{equation}
\mathbf{C}_\mathrm{1c}^{IG} = \int_{\Gamma} \mathbf{L}_{\mathrm{T}} \mathbf{R}_{\mathrm{1c}}^T \mathrm{d} \Gamma =  \int_{\Gamma} \mathbf{L}_{\mathrm{T}} \mathbf{L}_{\mathrm{1f}}^T \mathrm{d} \Gamma  \left(\mathbf{D}^{IGFE}_{\mathrm{1cf}}\right)^T = \mathbf{C}_\mathrm{1f}^{FE}  \left(\mathbf{D}^{IGFE}_{\mathrm{1cf}}\right)^T, \label{eq:C1cIG}
\end{equation}
with $\mathbf{D}^{IGFE}_{\mathrm{1cf}} =  \mathbf{D}^{IG}_{\mathrm{1cf}}  \mathbf{D}^{FE}_{\mathrm{1f}}$ and where $\mathbf{C}_\mathrm{1f}^{FE}$ is the mass matrix  associated with the refined global FE mesh of the FE interface. With these computations, Eqs.~\eqref{eq:nonint_global_dis}-\eqref{eq:nonint_local_dis} become:
\begin{equation}
\def\arraystretch{1.5}
\mathbf{K}_{\mathrm{1c}}^{IG} ~{\mathbf{u}_\mathrm{1c}^{IG}}^{(n)} = \mathbf{f}_{\mathrm{1c}}^{IG}  -\mathbf{D}^{IGFE}_{\mathrm{1cf}} { \mathbf{C}_\mathrm{1f}^{FE}}^T \boldsymbol{\lambda}^{(n-1)}+{\overline{\boldsymbol{\lambda}}_\mathrm{12c}^{IG}}^{(n-1)}\quad ;
\label{eq:nonint_global_dis_V1}
\end{equation}
\begin{equation}
\def\arraystretch{1.5}
\begin{bmatrix}
\mathbf{K}_{2}^{FE} & -{\mathbf{C}_2^{FE}}^T  \\
-\mathbf{C}_2^{FE} & \mathbf{0}  \\
\end{bmatrix}
\begin{pmatrix}
{\mathbf{u}_2^{FE}}^{(n)} \\ \boldsymbol{\lambda}^{(n)} 
\end{pmatrix}
=
\begin{pmatrix}
\mathbf{f}_2^{FE} \\  -\mathbf{C}_\mathrm{1f}^{FE}  \left(\mathbf{D}^{IGFE}_{\mathrm{1cf}}\right)^T {\mathbf{u}_{1}^{IG}}^{(n)}   \end{pmatrix},
\label{eq:nonint_local_dis_V1}
\end{equation}
which no longer exhibit operators that merge basis functions from IGA and FEM.

\begin{figure}[htb]{
		\centering
		\includegraphics[width=\linewidth]{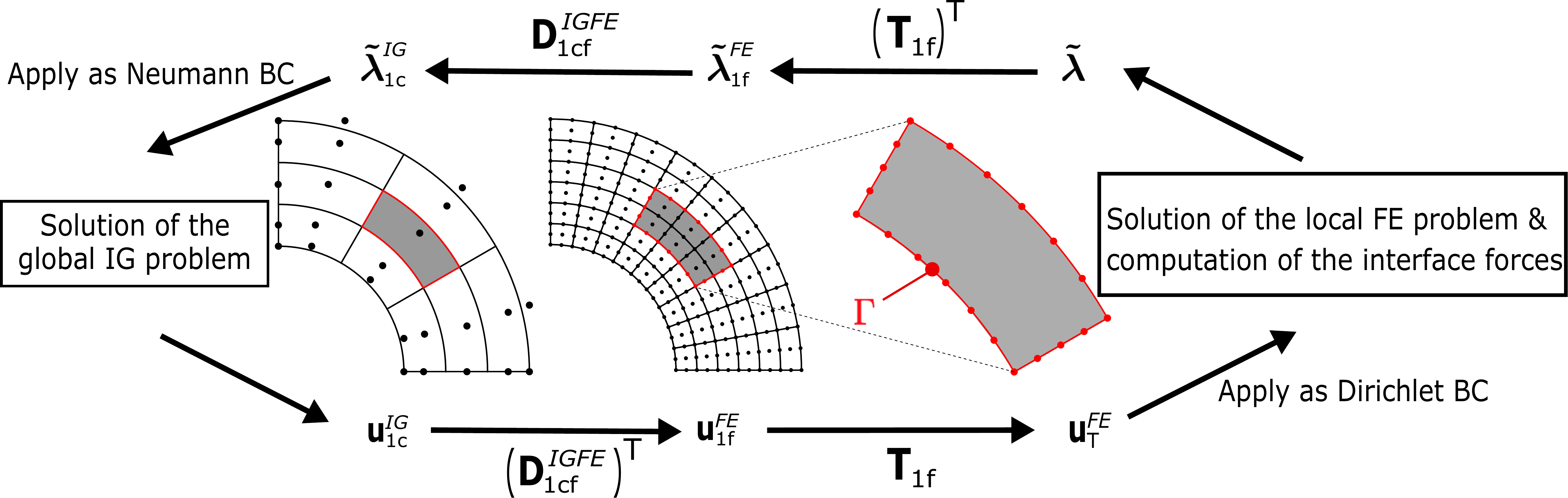}
		\caption{Communications between the global IG and local FE model through the fully non-invasive strategy.}
		\label{fig:comm_precise_global_local}}
\end{figure}

Now, the computation of the FE interface mass matrices $\mathbf{C}_\mathrm{1f}^{FE}$ and  $\mathbf{C}_2^{FE}$ may still not be straightforward in practice using commercial FE codes. However, as often performed in a transparent manner when coupling domains in FEM, these operators are not truly required here; the trace operators $\mathbf{T}_\mathrm{1f}$ and $\mathbf{T}_\mathrm{2}$ are actually sufficient. To highlight this, we first introduce the FE mass matrix of the interface:
\begin{equation}
\mathbf{C}_\mathrm{T}^{FE} = \int_{\Gamma} \mathbf{L}_{\mathrm{T}}  \mathbf{L}_{\mathrm{T}} ^T \mathrm{d} \Gamma,
\end{equation}
which is invertible (it symmetric definite positive by construction). Then, we take $\tilde{\boldsymbol{\lambda}}= \mathbf{C}_\mathrm{T}^{FE} \boldsymbol{\lambda}$ which has the dimension of a load vector, and we multiply the second part of Eq.~\eqref{eq:nonint_local_dis_V1} by $ \left( \mathbf{C}_\mathrm{T}^{FE} \right)^{-1}$. With these manipulations and making use of equalities:
\begin{equation}
\mathbf{T}_\mathrm{1f} = \left( \mathbf{C}_\mathrm{T}^{FE} \right)^{-1} \mathbf{C}_\mathrm{1f}^{FE} \quad ; \quad \mathbf{T}_\mathrm{2} = \left( \mathbf{C}_\mathrm{T}^{FE} \right)^{-1} \mathbf{C}_\mathrm{2}^{FE},
\end{equation}
Eqs.~\eqref{eq:nonint_global_dis_V1} and~\eqref{eq:nonint_local_dis_V1} read (modifications highlighted in color grey):
\begin{equation}
\def\arraystretch{1.5}
\mathbf{K}_{\mathrm{1c}}^{IG}~{\mathbf{u}_\mathrm{1c}^{IG}}^{(n)} =\mathbf{f}_{\mathrm{1c}}^{IG}  -\mathbf{D}^{IGFE}_{\mathrm{1cf}} \textcolor{black!40}{\mathbf{T}_{1f}^T ~\tilde{\boldsymbol{\lambda}}^{(n-1)}}+{\overline{\boldsymbol{\lambda}}_\mathrm{12c}^{IG}}^{(n-1)}\quad ;
\label{eq:nonint_global_dis_V2}
\end{equation}
\begin{equation}
\def\arraystretch{1.5}
\begin{bmatrix}
\mathbf{K}_{2}^{FE} & \textcolor{black!40}{-\mathbf{T}_{2}^T}  \\
\textcolor{black!40}{-\mathbf{T}_{2}} & \mathbf{0}  \\
\end{bmatrix}
\begin{pmatrix}
{\mathbf{u}_2^{FE}}^{(n)} \\ \textcolor{black!40}{\tilde{\boldsymbol{\lambda}}^{(n)} }
\end{pmatrix}
=
\begin{pmatrix}
\mathbf{f}_2^{FE} \\  \textcolor{black!40}{-\mathbf{T}_{1}}   \left(\mathbf{D}^{IGFE}_{\mathrm{1cf}}\right)^T {\mathbf{u}_{1}^{IG}}^{(n)}   \end{pmatrix}.
\label{eq:nonint_local_dis_V2}
\end{equation}
The trace operators being most of time available in FE codes, the coupling can now be carried out using only FE industrial packages, which is the reason why we characterize this method as fully non-invasive. For completeness, Fig.~\ref{fig:comm_precise_global_local} further illustrates the communications between the IG global and local FE model in line with Eqs.~\eqref{eq:nonint_global_dis_V2}-\eqref{eq:nonint_local_dis_V2}.

\begin{remark}
	Pushing forward the reasoning by applying Eq.~\eqref{eq:reduced_order}, an implementation using only FE industrial packages for both the coupling and classical stiffness and load vector operators can also be proposed if the users only have at their disposal FE codes. Indeed, we can compute:
	\begin{subequations}
		\begin{align}
		&\mathbf{K}_{\mathrm{1c}}^{IG} =  \mathbf{D}^{IGFE}_{\mathrm{1cf}}\mathbf{K}_{\mathrm{1f}}^{FE}\left(\mathbf{D}^{IGFE}_{\mathrm{1cf}}\right)^T \quad \mathrm{and} \quad \mathbf{f}_{\mathrm{1c}}^{IG} =  \mathbf{D}^{IGFE}_{\mathrm{1cf}}\mathbf{f}_{\mathrm{1f}}^{FE} ; \label{eq:K1cIG} \\
		&\overline{\boldsymbol{\lambda}}_\mathrm{12c}^{IG} =  \mathbf{D}^{IGFE}_{\mathrm{1cf}} \overline{\boldsymbol{\lambda}}_\mathrm{12f}^{FE} \quad \mathrm{with} \quad \overline{\boldsymbol{\lambda}}_\mathrm{12f}^{FE} = \overline{\mathbf{K}}_\mathrm{12f}^{FE} \mathbf{u}_\mathrm{1f}^{FE}- \overline{\mathbf{f}}_\mathrm{12f}^{FE}, \quad  \mathbf{u}_\mathrm{1f}^{FE} = \left(\mathbf{D}^{IGFE}_{\mathrm{1cf}}\right)^T \mathbf{u}_\mathrm{1c}^{IG} ; \label{eq:r12cIG}
		\end{align}
	\end{subequations}
	which yields the iterative following process:
	\begin{equation}
	\def\arraystretch{1.5}
	\mathbf{D}^{IGFE}_{\mathrm{1cf}}\mathbf{K}_{\mathrm{1f}}^{FE}\left(\mathbf{D}^{IGFE}_{\mathrm{1cf}}\right)^T {\mathbf{u}_\mathrm{1c}^{IG}}^{(n)}  =\mathbf{D}^{IGFE}_{\mathrm{1cf}}\mathbf{f}_{\mathrm{1f}}^{FE}  -\mathbf{D}^{IGFE}_{\mathrm{1cf}} \mathbf{T}_{1f}^T ~\tilde{\boldsymbol{\lambda}}^{(n-1)}+\mathbf{D}^{IGFE}_{\mathrm{1cf}} {\overline{\boldsymbol{\lambda}}_\mathrm{12f}^{FE}}^{(n-1)};
	\label{eq:nonint_global_dis_V3}
	\end{equation}
	\begin{equation}
	\def\arraystretch{1.5}
	\begin{bmatrix}
	\mathbf{K}_{2}^{FE} & -\mathbf{T}_{2}^T \\
	-\mathbf{T}_{2} & \mathbf{0}  \\
	\end{bmatrix}
	\begin{pmatrix}
	{\mathbf{u}_2^{FE}}^{(n)} \\ \tilde{\boldsymbol{\lambda}}^{(n)}
	\end{pmatrix}
	=
	\begin{pmatrix}
	\mathbf{f}_2^{FE} \\  -\mathbf{T}_{1}  \left(\mathbf{D}^{IGFE}_{\mathrm{1cf}}\right)^T {\mathbf{u}_{1}^{IG}}^{(n)}   \end{pmatrix},
	\label{eq:nonint_local_dis_V3}
	\end{equation}
	where only FE operators are required (see Fig.~\ref{fig:gIGAlFEM_papier_principe} for the notations).
	\label{rk:FEM/FEM}
\end{remark}

\begin{remark}
	With the proposed method, it is possible to take a higher-order IG model (\emph{i.e.}, $p>1$) with a classic low-order FE mesh (\emph{i.e.}, $p=1$). In this case, we take as many nodes at the interface for the refined global FE model (see Fig.~\ref{fig:gIGAlFEM_papier_principe}(c)) as for the local FE model (see  Fig.~\ref{fig:gIGAlFEM_papier_principe}(e)). This means that, for instance, for quadratic IGA versus linear FEM, there are two FE elements in front of one IG element along the interface. Two approximations are thus performed implicitly: (i) the (possibly curved) interface of the IG model is facetted between each of the interface nodes for the local model, (ii) $\left( \mathbf{C}_\mathrm{T}^{FE} \right)^{-1} \mathbf{C}_\mathrm{1f}^{FE}$ is not strictly equal to $\mathbf{T}_\mathrm{1f}$ any more (the mass matrix of a quadratic FE interface is approximated by the mass matrix of a linear FE interface with twice more elements). However, our results seem to indicate that the error related to such approximations is very low compared to that associated with the discretization (see section~\ref{sub:res_lin_2D} and Fig.~\ref{fig:lin:lin:raf}).
	\label{rk:Global-quad/Local-lin}
\end{remark}



\section{Numerical results}
\label{sec:results}

To assess the performance of the developed non-invasive hybrid global-IGA/local-FEM algorithm, we now present a series of numerical experiments that cover 2D and 3D simulations with different local behaviors, such as cracks, contact and delamination. All the implementations have been carried out using the open-source FE industrial software package Code\_Aster~\cite{aster} developed by the EDF R\&D company. No IG codes have been used; we consider the case where we have at our disposal only the FE code Code\_Aster. We thus more precisely implemented algorithm~\eqref{eq:nonint_global_dis_V3}-\eqref{eq:nonint_local_dis_V3} and limited ourselves to quadratic spline functions since Code\_Aster does not go beyond second-order Lagrange finite elements. Yet, we underline that the proposed implementation schemes~\eqref{eq:nonint_global_dis_V1}-\eqref{eq:nonint_local_dis_V1} and~\eqref{eq:nonint_global_dis_V2}-\eqref{eq:nonint_local_dis_V2} could be applied to higher-order splines if one has an IG code in hand to compute the IG stiffness and load vector operators. The automatic procedure described in Fig.~\ref{fig:gIGAlFEM_papier_principe} for the construction of the conforming global/local discretization was performed using the mesh generator Salome-Meca~\cite{salome} included in Code\_Aster. Finally, every computed fields (displacement, stress) are expressed in terms of FE quantities (see Eq.~\eqref{eq:link_disp}) so that the standard Code\_Aster post-processing functionalities are used to visualize the results. In the illustrations, we keep the notations introduced previously in the paper; in particular, domain $\Omega_1 = \Omega_{11} \cup \Omega_{12} \cup \Gamma$ characterizes the global IG model, and the local FE model of domain $\Omega_2$ is expected to replace the global IG model in sub-domain $\Omega_{12}$.

\subsection{Linear elastic 2D curved beam}
\label{sub:res_lin_2D}

\begin{figure}[ht!]
	\centering
	\subfloat[Problem description and discretization. \label{fig:lin:casTest}]{\includegraphics[width=0.4\textwidth]{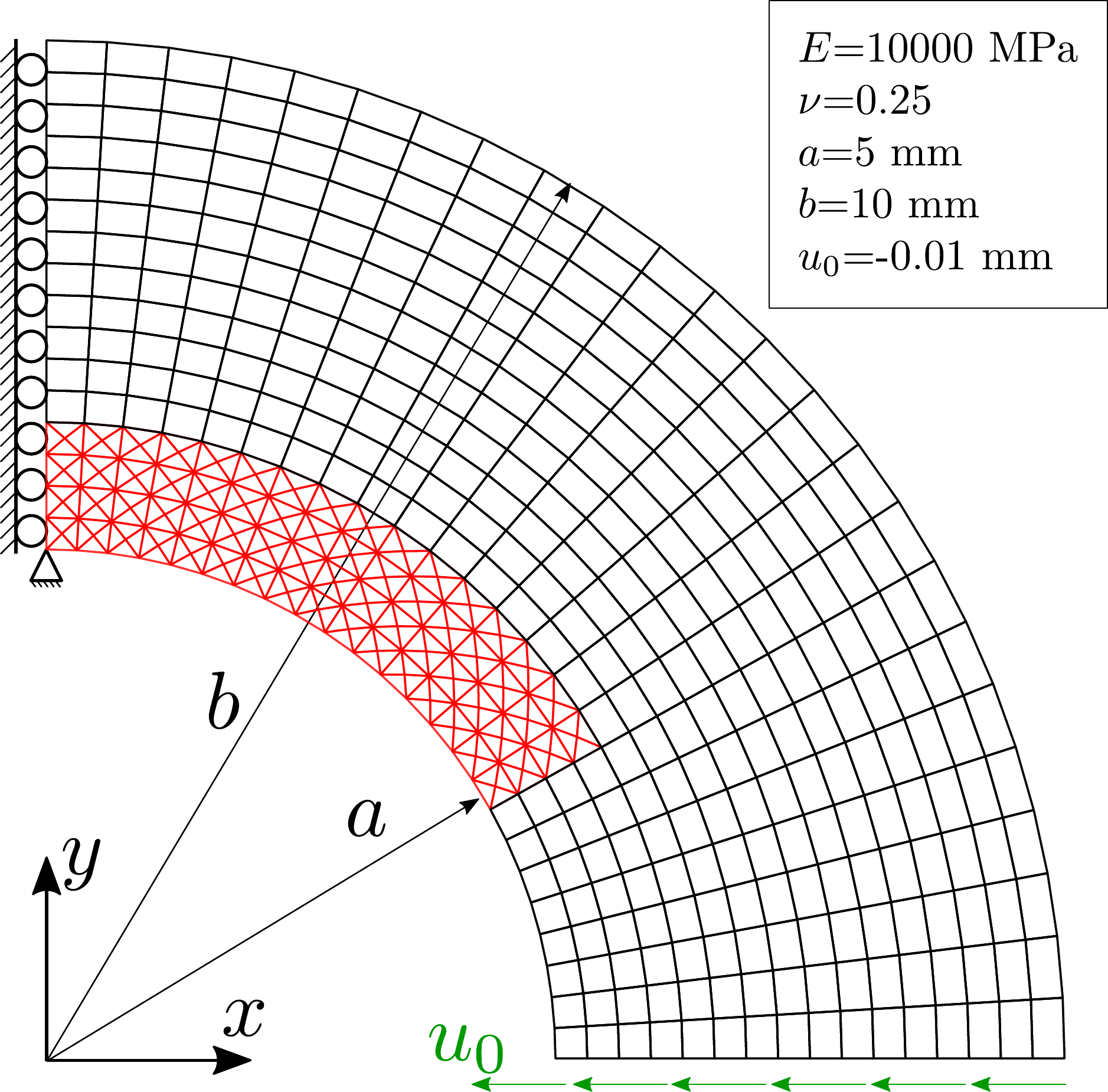}}\hfill
	\subfloat[Displacement field.\label{fig:lin:disp}]{\includegraphics[width=0.38\textwidth]{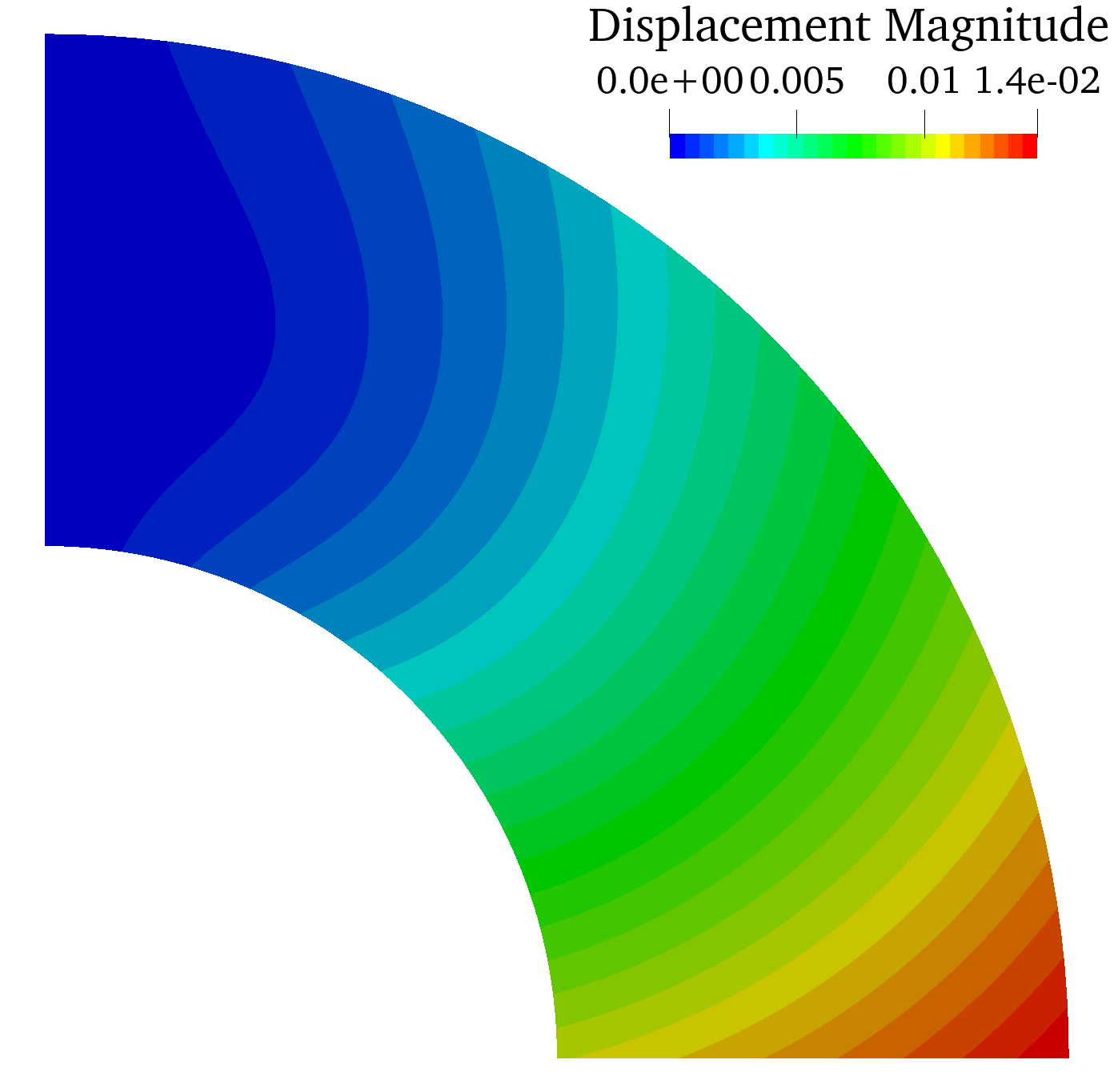}}\hfill
	\subfloat[Von Mises stress.\label{fig:lin:vm}]{\includegraphics[width=0.38\textwidth]{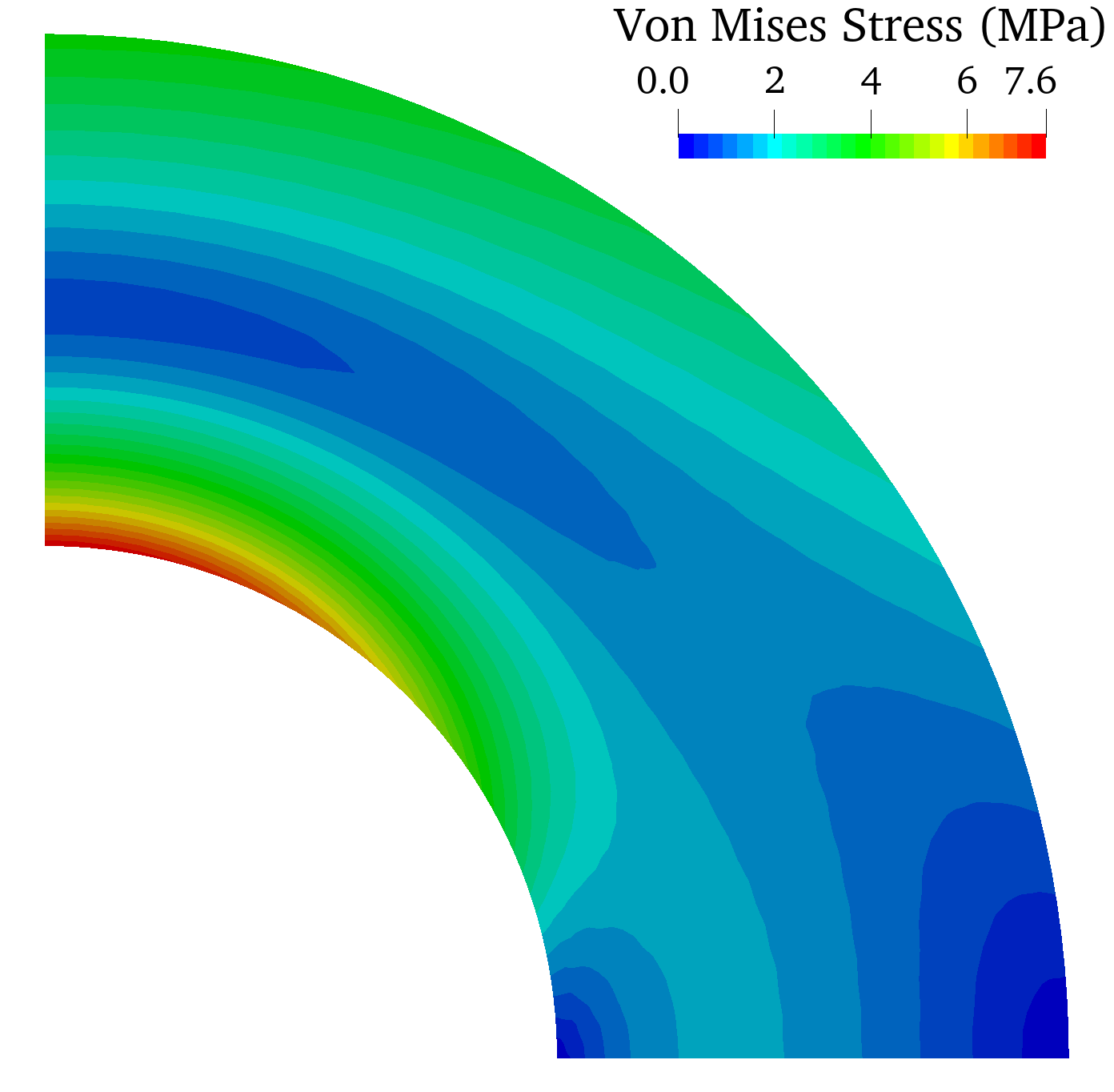}}
	\caption{Global/local non-invasive analysis of the linear curved beam problem (NURBS mesh composed of quadratic 24 (circumferential direction) $\times$ 16 (radial direction) elements for $\Omega_1$, and standard FE mesh composed of 128 T6 elements for $\Omega_2$.}
\end{figure}

The first example consists of a 2D linear elastic curved beam subjected to end shear adapted from~\cite{zienck_book}. Such an example has been widely used in IGA to assess the performance a method. The global geometry was perfectly generated using a single NURBS patch composed of only one quadratic element. The problem, together with the proposed global/local discretization, is illustrated in Fig.~\ref{fig:lin:casTest}. The plane stress assumption was performed and a constant horizontal displacement of $u_0= -0.01$ mm was prescribed over the lower beam boundary. In a small part of the bottom-left corner of the structure, where stress concentrations may appear, the global quadratic NURBS model was meant to be substituted by a local standard FE model composed of quadratic triangles (\emph{i.e.}, T6 triangles). More precisely, each initial quadratic NURBS element in $\Omega_{12}$ was replaced by 4 T6 elements in $\Omega_2$. The results obtained by performing algorithm~\eqref{eq:nonint_global_dis_V3}-\eqref{eq:nonint_local_dis_V3}  with the discretization of Fig.~\ref{fig:lin:casTest} are shown in Figs.~\ref{fig:lin:disp} and~\ref{fig:lin:vm} in terms of displacement and of Von Mises stress, respectively. We note that it is the converged solution in $\Omega_{11}\cup\Omega_2$ that is mapped (the fictitious prolongation of the global solution over $\Omega_{12}$ is not represented). For all examples, we will perform this way for the visualization. On this simple case, the iterative non-invasive algorithm converges very quickly without acceleration techniques: 3 iterations were needed with a stopping criterion based on the equilibrium of the interface global and local reaction forces (tolerance of $10^{-8}$ here). The solution appears smooth and in a good agreement with~\cite{zienck_book}.

\begin{figure}[ht!]
	\centering
	\includegraphics[scale=0.9]{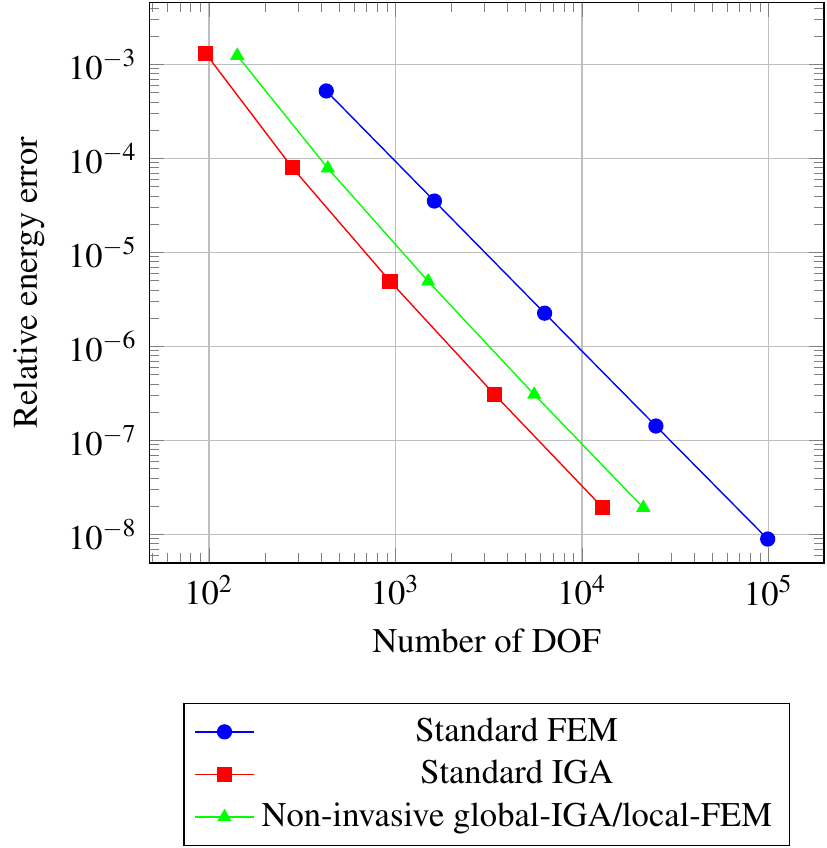}
	\caption{Convergence of the relative energy error for the linear elastic 2D circular beam (quadratic case).}
	\label{fig:fig:lin:convQuad}
\end{figure}

To go further, the convergence of the method with the refinement of the mesh was studied. In order to do so, the computation depicted in Fig.~\ref{fig:lin:casTest} was repeated for several global/local discretizations. Starting with a global quadratic NURBS mesh of $6$ (circumferential direction) $\times$ $4$ (radial direction) elements, the refinement was increased to reach $12\times 8$, $24 \times 16$, $48 \times 32$ and $96 \times 64$ elements. We kept the same region $\Omega_{12}$ and the replacement of 1 NURBS element by 4 T6 triangles for the local model for each discretization. We proceeded in the same way as in~\cite{zienck_book}; that is, the convergence behavior of the strain energy is considered through to computation of the relative energy error:
%
%
\begin{equation}
\mathrm{Err}^h = \frac{\lvert E^{\mathrm{ref}}-E^h \rvert}{E^{\mathrm{ref}}}
\end{equation}
where $E^\mathrm{ref}$ denotes the reference exact strain energy and $E^h$ the strain energy of the discrete model. The convergence curve is given versus the number of DOF in Fig.~\ref{fig:fig:lin:convQuad} (see green curve). The number of DOF was computed as the sum of the global IG DOF and of the local FE DOF. For comparison purpose, the convergence curves of the equivalent single-model solutions are also plotted: "Standard IGA" represents the solution when considering the global model everywhere and "Standard FEM" corresponds to the solution when all the NURBS elements are replaced by 4 T6 triangles. The results show that the same rate of convergence was achieved with the proposed hybrid IGA/FEM scheme as with the reference solutions, which accounts for the accuracy of our method. More specifically, it can be noticed that for a given mesh refinement, the errors are about the same for the three solutions. Only the number of DOF changes: it decreases when IGA is used. This illustrates the increased per-DOF accuracy of IGA and is totally consistent with our interpretation on IGA as a projection of FEM onto a reduced, regular basis. The solution here being smooth, it is well captured with IGA as with FEM, but IGA comes with less DOF due to its higher regularity. Our hybrid IGA/FEM solution obviously appears between the two reference curves since the DOF of the global and local models are summed. Let us underline here that our way of counting the DOF has not a concrete meaning from a computational cost point of view since the IG and FE problems are solved separately in our non-invasive strategy. One could have chosen to take the maximum number of DOF between the global and local models, \emph{i.e.} the same number of DOF as for the standard IGA solution which would have led to the superposition of the "Standard IGA" and "Non-invasive global-IGA/local-FEM" curves.

\begin{figure}[ht!]
	\centering
	\subfloat[Local T3 triangles (in black) \label{fig:lin:lin:mesh} compared to the initial quadratic IG mesh over region $\Omega_{12}$ (in red).]{\includegraphics[width=0.3\textwidth]{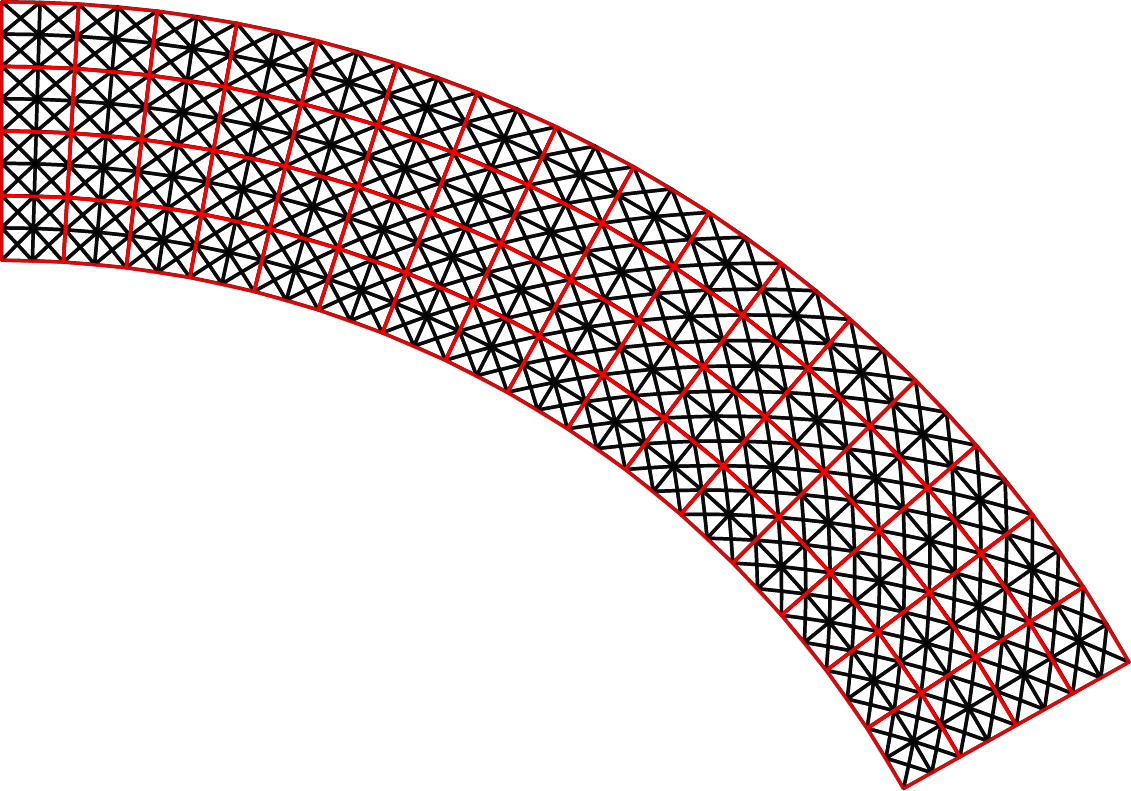}}\hfill
	\subfloat[Displacement field.\label{fig:lin:lin:displ}]{\includegraphics[width=0.38\textwidth]{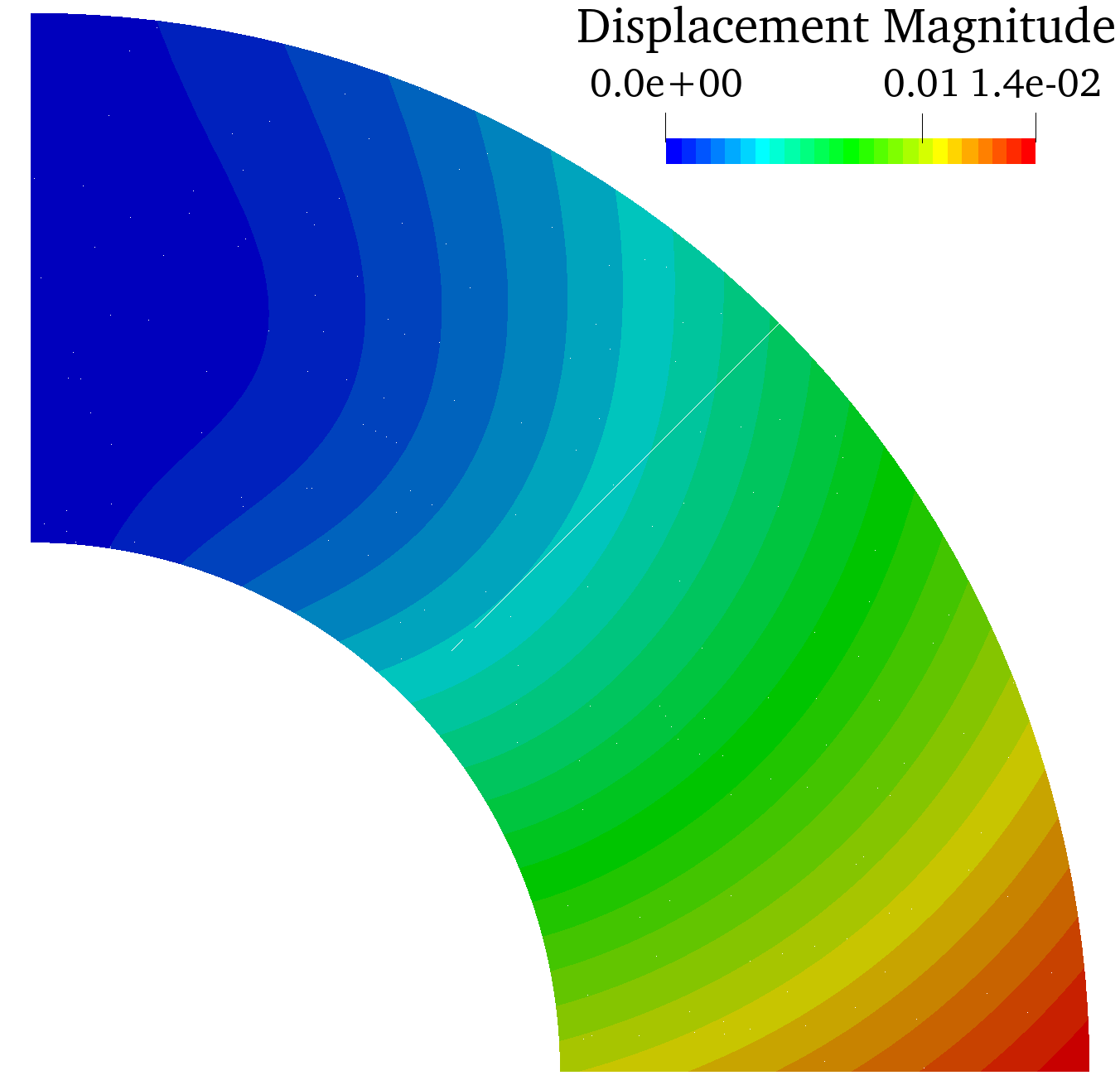}}\hfill
	\subfloat[Von Mises stress.\label{fig:lin:lin:vm}]{\includegraphics[width=0.38\textwidth]{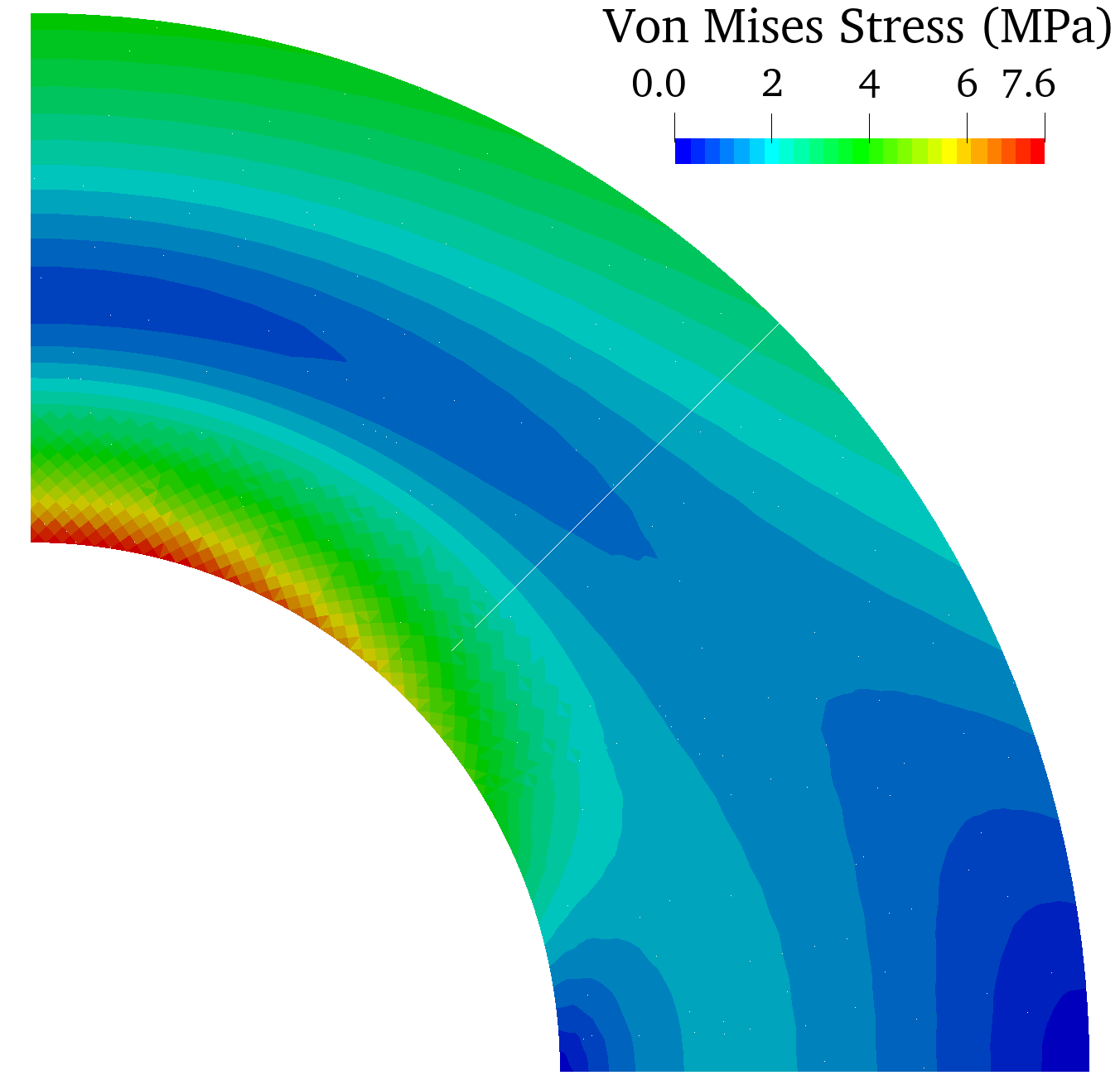}}
	\caption{Global/local non-invasive analysis of the linear curved beam problem (NURBS mesh composed of quadratic 24 (circumferential direction) $\times$ 12 (radial direction) elements for $\Omega_1$, and standard FE mesh composed of 512 T3 elements for $\Omega_2$.}
	\label{fig:my_label}
\end{figure}

\begin{figure}
	\centering
	\includegraphics[scale=0.9]{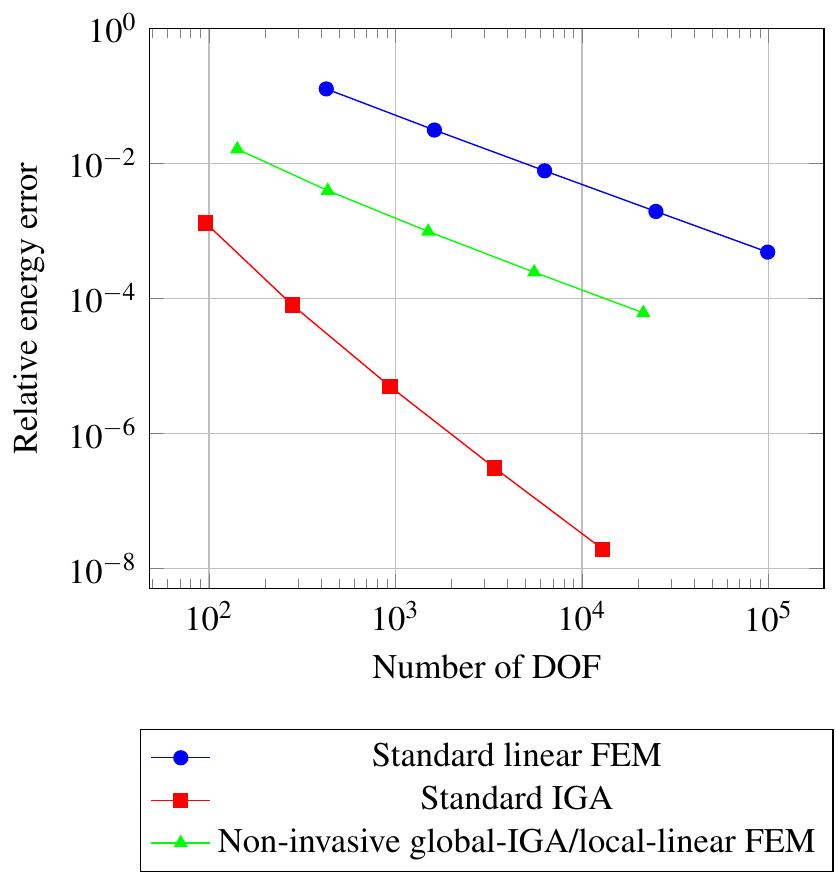}
	\caption{Convergence of the relative energy error  for the linear elastic 2D circular beam with global quadratic NURBS elements and local T3 triangles.}
	\label{fig:lin:lin:raf}
\end{figure}

For completeness on this test case, we finally carried out the same numerical experiments but with a local FE model composed of linear triangles (\emph{i.e.}, T3 triangles). For this purpose, the strategy described in remark~\ref{rk:Global-quad/Local-lin} was applied: two FE elements were put in front of one NURBS element along the interface. Fig.~\ref{fig:lin:lin:mesh} shows a zoom on the local region when considering a global model made of $24 \times 16$ elements. This time, 16 T3 triangles replaced 1 NURBS element in the local region. The results for this discretization are given in Fig.~\ref{fig:lin:lin:displ} and~\ref{fig:lin:lin:vm} in terms of displacement and of Von Mises stress, respectively. Of course, some discontinuities for the stress can be observed in the local region since this field is now piecewise constant in this area, but the solution still appears in a good agreement with~\cite{zienck_book}. The convergence curve in terms of relative energy error is then plotted in Fig.~\ref{fig:lin:lin:raf} along with the equivalent single-model solutions; that are, the solution when considering the global quadratic IG model everywhere and the solution when all the NURBS elements are replaced by 16 T3 triangles. As expected, the convergence rate of our hybrid quadratic-IGA/linear-FEM scheme is now driven by the linear-FEM part of the solution (same convergence rate for the coupled solution as for the linear-FEM case). However, the coupled solution appears much more accurate than the standard linear-FEM one (drastic reduction of the constant factor) which is due to the higher order and higher regularity of the global model. This numerical experiment validates the proposed procedure to couple higher-order IGA with standard linear FEM (see again remark~\ref{rk:Global-quad/Local-lin}) and further confirms the interest of making use of IGA for the global response to reach an increased per-DOF accuracy.

\subsection{2D curved beam with holes, cracks and contact}
\label{sub:res_nonlin_2D}

\begin{figure}[ht!]
	\centering
	\includegraphics[width=0.9\textwidth]{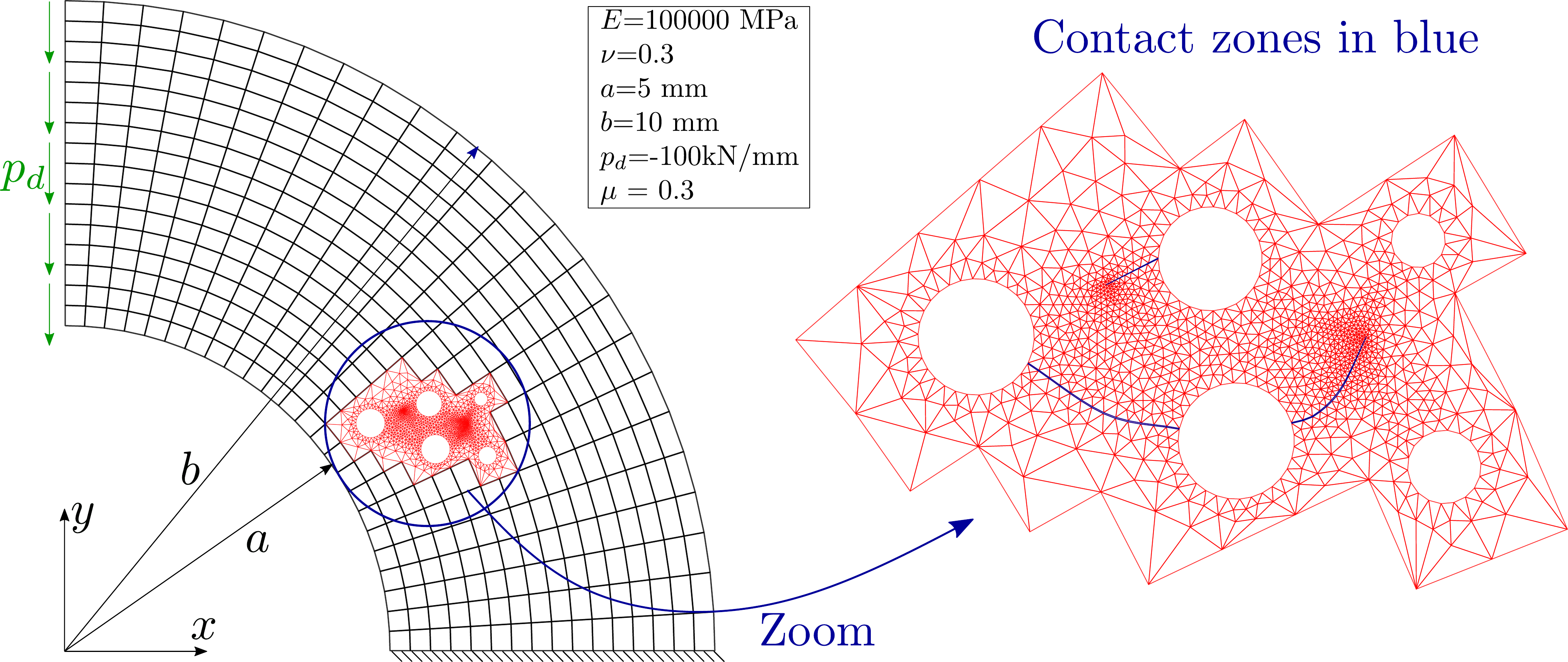}
	\caption{Non-invasive introduction of holes, cracks and frictional contact between the lips of the cracks in an initial 2D NURBS beam with the proposed hybrid global-IGA/local-FEM methodology (the cracks, where a contact model is applied, are colored in blue).}
	\label{fig:cas-2D-contact}
\end{figure}

With the second example, we illustrate the potential of our methodology to include geometrical details along with non-linear behaviors within a global NURBS model. More precisely, the global linear elastic 2D curved beam model of the previous test case was recycled. This time, it was clamped on its lower part and subject to a vertical distributed load on its left edge, see Fig.~\ref{fig:cas-2D-contact}(left). The global mesh was composed of 24 $\times$ 16 quadratic NURBS elements. This model was locally enriched by a specific FE mesh made of quadratic triangles (\emph{i.e.}, T6 triangles) incorporating several holes and cracks. In addition, frictional contact was modeled between the lips of the cracks. The cracks were initially open (initial gap between the crack lips of about 0.004 mm). A zoom on the local model is performed in Fig.~\ref{fig:cas-2D-contact}(right). Let us underline that solving such a problem within the sole IG framework would be delicate due to the geometric complexity of the local region (perhaps, some advanced IG immersed technologies~\cite{Ruess14,Wei21,Wang21} or extended~\cite{Luycker11,Yuan21,Fathi21} or phase-field procedures~\cite{Borden14,Proserpio20,Paul20} should be required). On the contrary, realizing a locally boundary-fitted mesh is straightforward with standard FEM. Besides, the displacement has to be discontinuous on either side of the cracks which is natural with a FE mesh that fits the cracks.

\begin{figure}[ht!]
	\centering
	\subfloat[Von Mises stress.\label{fig:vonmises-contact}]{\includegraphics[width=0.72\textwidth]{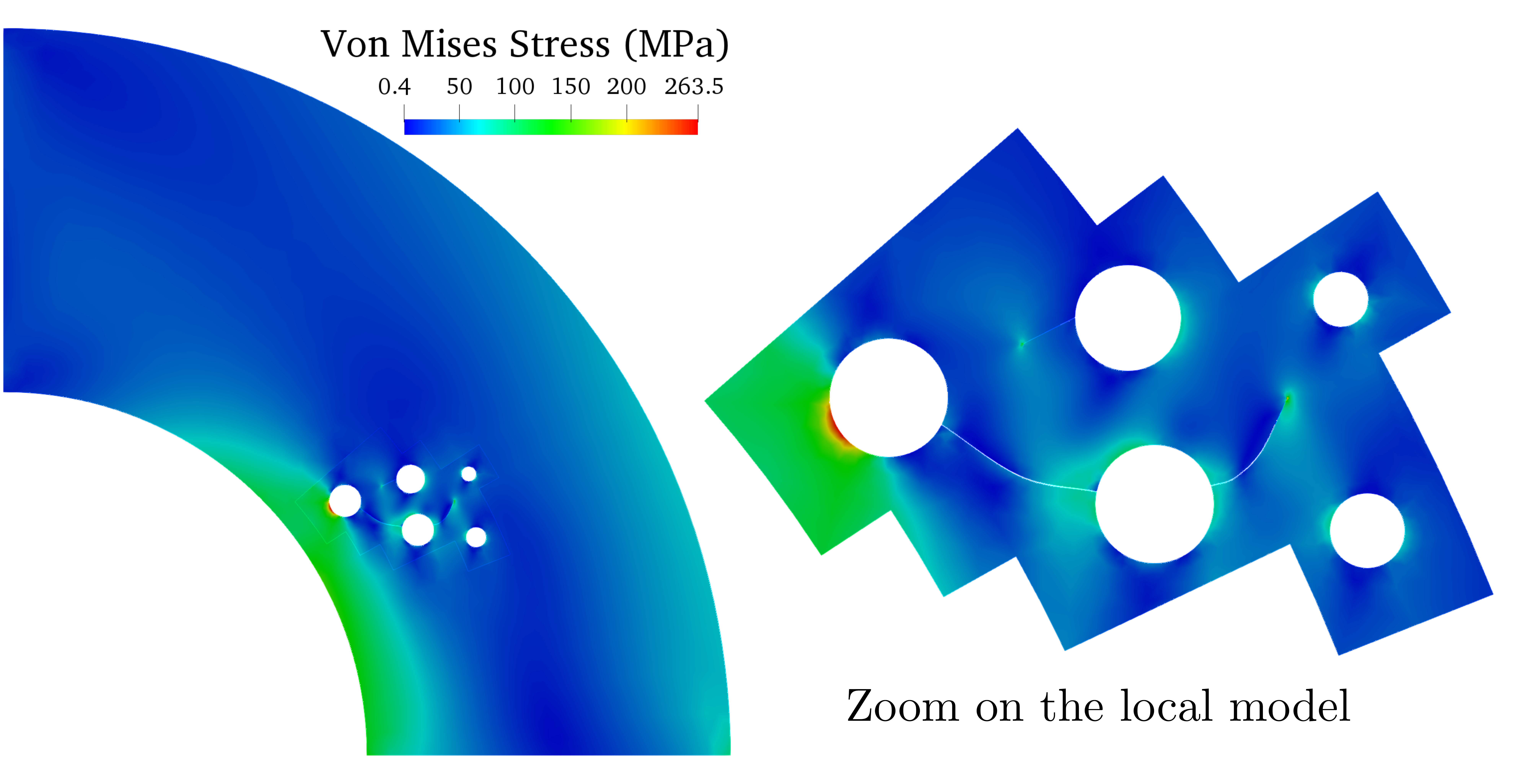}} \qquad
	\subfloat[Zoom on contact zones of the deformed configuration (scale factor 1) - contact is reached on the red zones.\label{fig:contact-zones-2D}]{\includegraphics[width=0.37\textwidth]{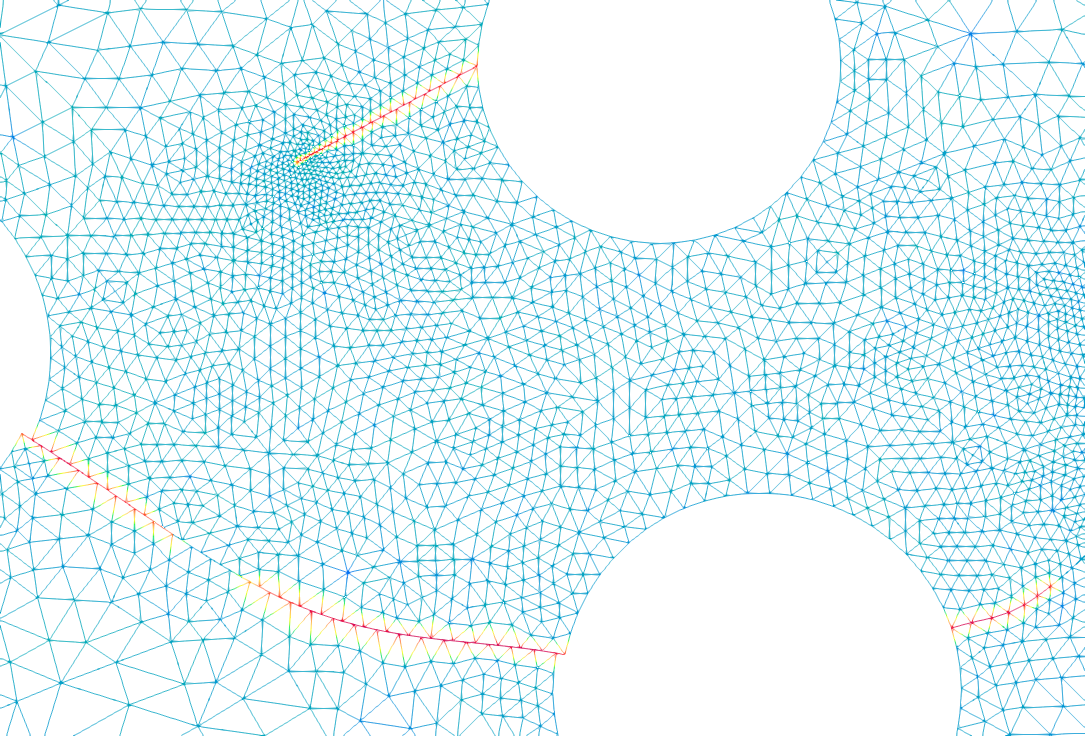}} \qquad
	\subfloat[Relative displacement discrepancy (\%) with respect to an equivalent full FE solution. \label{fig:diff-displ}]{\includegraphics[width=0.37\textwidth]{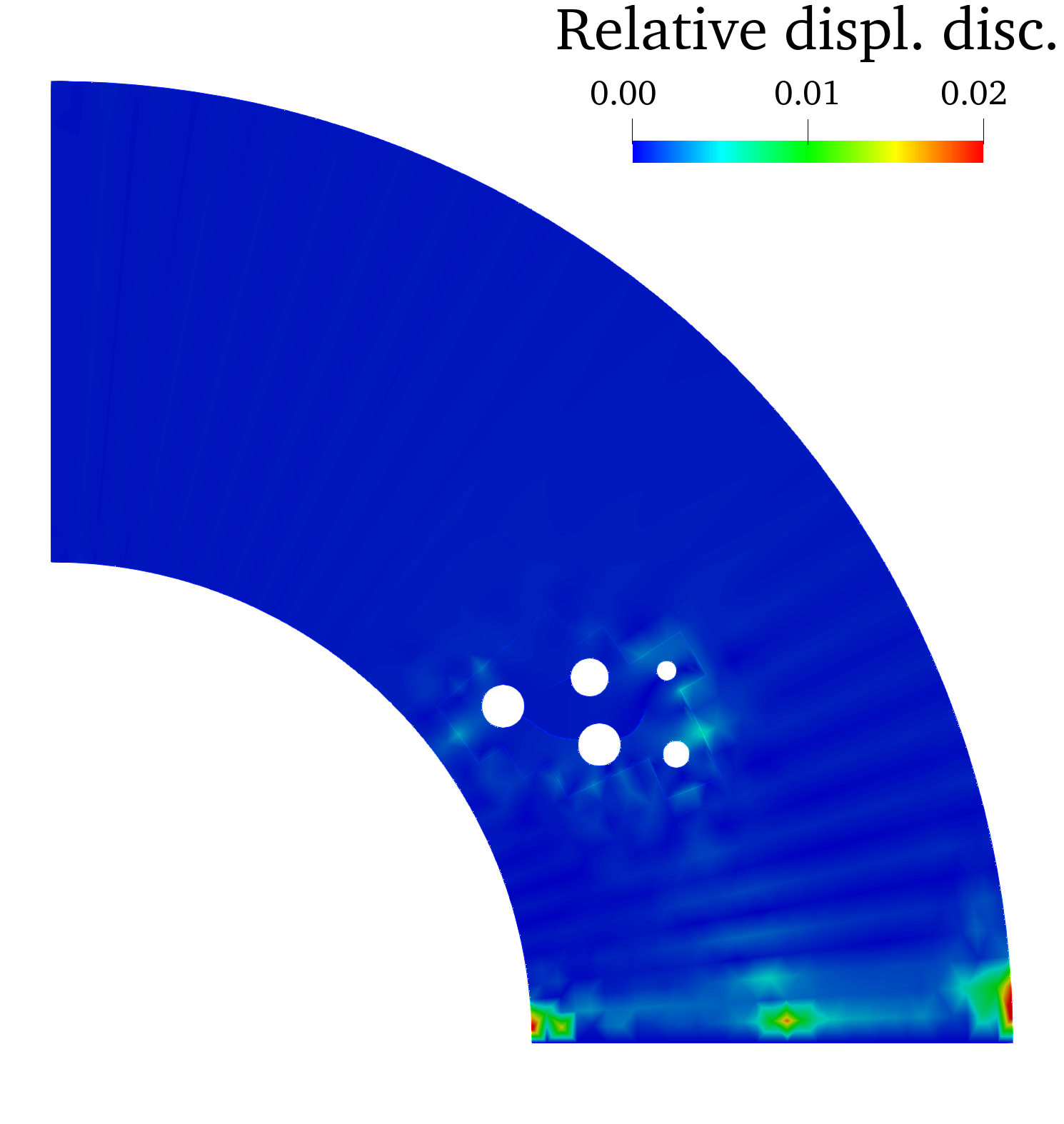}}
	\caption{Solution obtained for the 2D curved beam with holes, cracks and contact with the proposed hybrid global-IGA/local-FEM algorithm.}
	\label{fig:res-contact}
\end{figure}

The Von Mises stress map obtained once the non-invasive algorithm has converged is given in Fig.~\ref{fig:vonmises-contact}. The transition of the stress at the global/local interface appears smooth, although the global/local discretization is $C^0$ at that location, which confirms the accuracy of our hybrid coupling. Stress concentrations are observed close to the holes and at the crack tips which is mechanically sound. Furthermore, the contact zones on the crack lips are highlighted in Fig.~\ref{fig:contact-zones-2D}. The cracks are closing which seems to be consistent with the applied load and Dirichlet boundary conditions. More precisely, the whole top crack closes while only subparts of the other cracks are in contact. In addition, slight sliding can be observed on the large bottom crack. For completeness, the global-IGA/local-FEM displacement obtained through our non-invasive strategy was compared to the displacement field obtained with a full FE discretization of the whole problem (the FE mesh that allows to recover the IG global solution by projection (see Eq.~\eqref{eq:K1cIG}) was used in $\Omega_{11}$). Since the solutions are very close, the relative discrepancy between the hybrid solution and the FE one is plotted in Fig.~\ref{fig:diff-displ}. Less than 2\% of local mismatch can be observed while the two solutions come from different approximation subspaces in $\Omega_{11}$. This result confirms that IGA is sufficient to accurately capture the global response, even when the local one exhibits some discontinuities.


\begin{figure}
	\centering
	\includegraphics[scale=0.9]{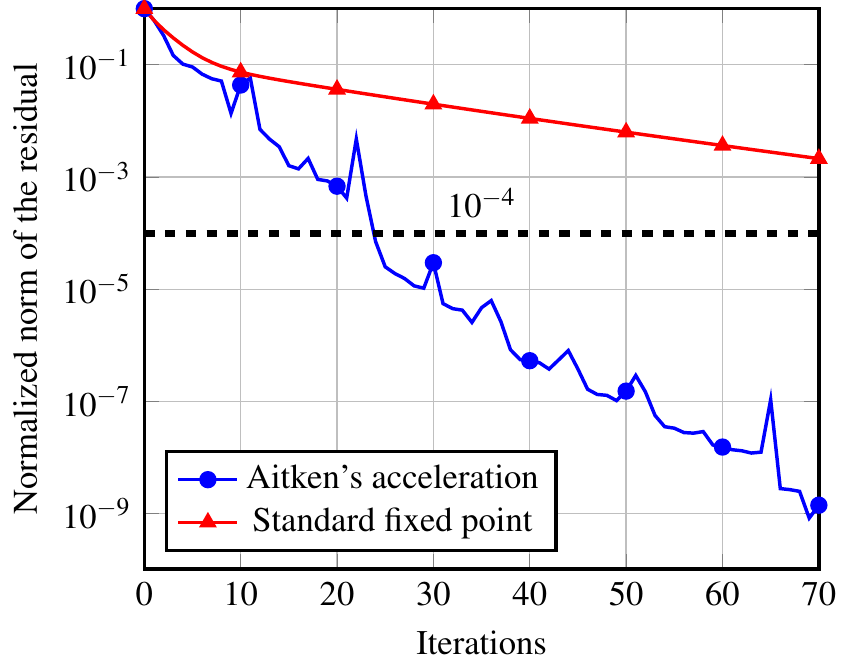}
	\caption{Convergence of the non-invasive global/local algorithm for the 2D curved beam with holes, cracks and contact.}
	\label{fig:conv-contact}
\end{figure}

Finally, the convergence of the non-invasive algorithm was investigated in Fig.~\ref{fig:conv-contact}. For this example, the Aitken's dynamic relaxation acceleration ~\cite{Duval16,Gosselet18} seems to be necessary, which was expected since the stiffness gap between domains $\Omega_{12}$ and $\Omega_2$ was significant. Taking advantage of the Aitken's update, we were able to make the number of iterations relatively low: a residual of $10^{-4}$ was obtained in 23 iterations.

\subsection{2D plate with multiple inclusions and delamination}
\label{sub:2D_inclusions_delamination}

The third example constitutes a 2D illustration of a composite material made of a matrix and several inclusions inside. The objective was to show that our approach is able to deal with cohesive zones in the local model in order compute the delamination at the inclusion-to-matrix interfaces. Moreover, this test case will allow to highlight another attractive property of our non-invasive algorithm: it results in an efficient non-linear domain decomposition solver when several local models are considered.

\begin{figure}[ht!]
	\centering
	\includegraphics[width=0.85\textwidth]{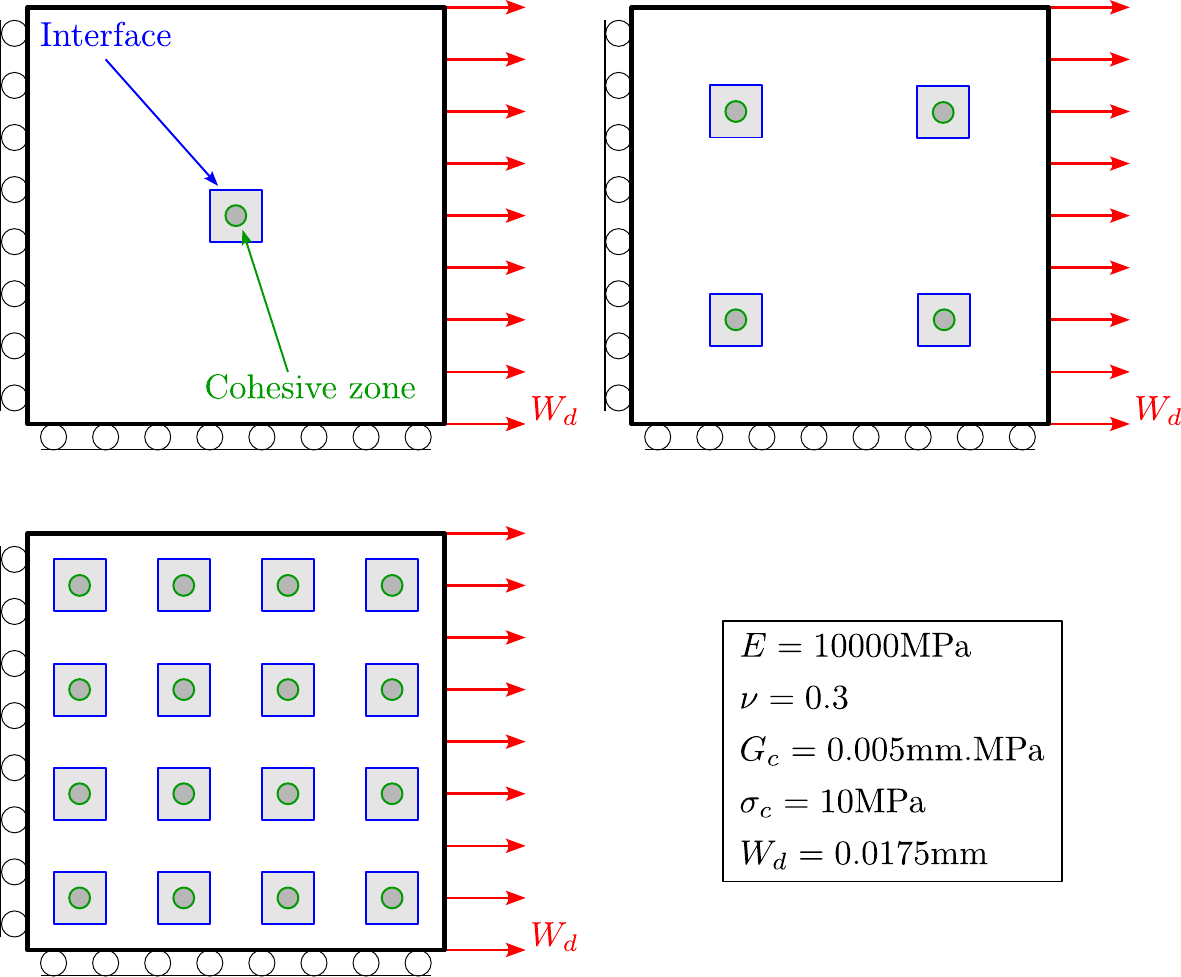}
	\caption{Description of the 2D plate problem with multiple inclusions and delamination. 1, 4 or 16 inclusions are considered. The global/local interfaces are underlined in blue while the cohesive interfaces are plotted in green. One local FE model includes one inclusion plus a matrix region surrounding the inclusion, and thus one cohesive interface. The different model parameters are indicated on the right.}
	\label{fig:caseInclusion}
\end{figure}

\begin{figure}[ht!]
	\centering
	\includegraphics[width=0.5\textwidth]{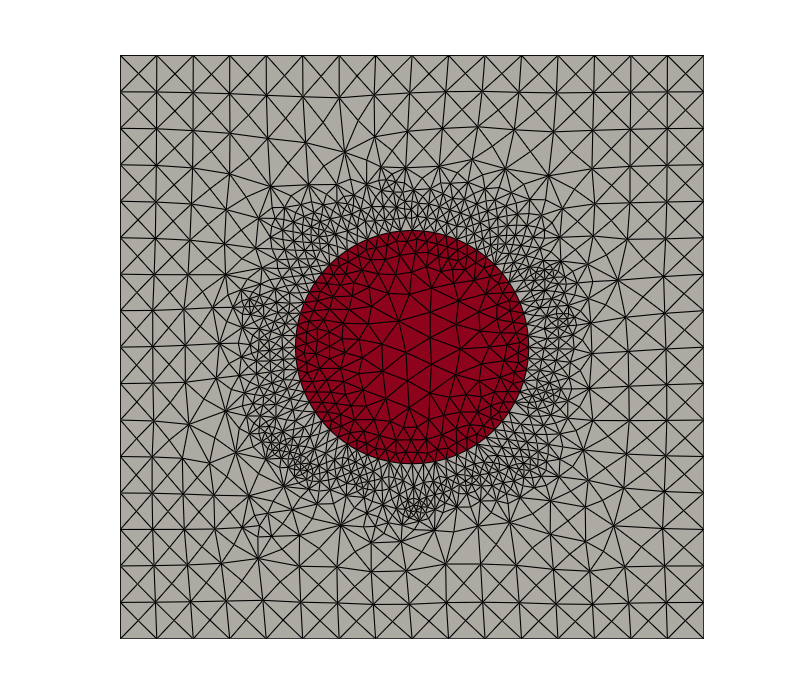}
	\caption{Local FE meshes considered for the plate with multiple inclusions problem. It is made of $T3$ triangles. The inclusion is in red and cohesive elements are incorporated at the inclusion-to-matrix interface. On the boundary of the complete local model, \emph{i.e.} at the global/local interface, 2 $T3$ triangles are put in front of one element of the intermediary refined global model, following remark~\ref{rk:Global-quad/Local-lin}.}
	\label{fig:localMeshInclusion}
\end{figure}

\begin{figure}[ht!]
	\centering
	\includegraphics[width=0.45\textwidth]{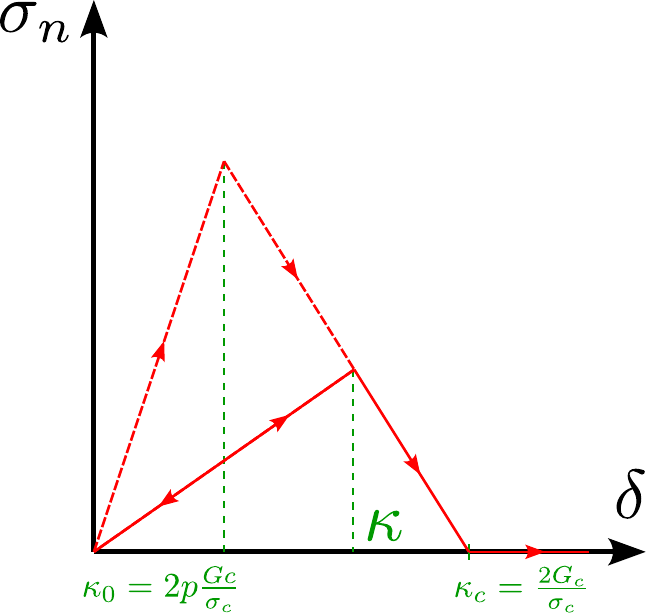}
	\caption{Bilinear law for the cohesive elements at the inclusion-to-matrix interfaces.}
	\label{fig:bilineal}
\end{figure}

More precisely, the structure was composed of a square shape matrix containing 1, 4 or 16 (uniformly distributed) circular inclusions, as depicted in Fig.~\ref{fig:caseInclusion}. Symmetry boundary conditions were applied on the left and bottom edge of the structure, while a constant horizontal displacement was prescribed over the right edge, which results in a structure globally subjected to traction. The global model, which obviously describe here the square shape matrix (without inclusions), was discretized with $16\times16$ quadratic IG elements. This IG model was then enriched by as many local FE models as there are inclusions. The associated local FE meshes were the same for all the inclusions and they can be seen on Fig.~\ref{fig:localMeshInclusion}. They all included one inclusion along with a matrix region surrounding the inclusion. This allows to incorporate a cohesive zone at the inclusion-to-matrix interface within the local FE model. To obtain conforming interfaces between the global model and the local ones, we made use of the strategy depicted in Fig.~\ref{fig:gIGAlFEM_papier_principe} with, this time, a spline refinement operator $\mathbf{D}^{IG}_{\mathrm{1cf}}$. In addition, we followed remark~\ref{rk:Global-quad/Local-lin} to adopt a global quadratic-IGA/local linear-FE modeling. We first refined the IG global model 4 times and applied the FEM-to-IGA bridge to obtain a refined global FE model made of 64$\times$64 quadratic FE elements. We then extracted the FE nodes on the global/local interface and built the FE mesh of Fig.~\ref{fig:localMeshInclusion} by putting $T3$ elements between each interface nodes. We eventually considered the same material for the fibers and the matrix: a linear elastic isotropic material with the Young modulus $E=10000$ MPa and Poisson ratio $\nu=0.3$. The cohesive elements at the inclusion-to-matrix interfaces followed a standard bilinear law, see Fig.~\ref{fig:bilineal}. This law is described by 3 parameters : $\sigma_c$ the critical stress, $G_c$ the density of the critical energy of the material and $p$ a penalisation coefficient to control the stiffness of the undamaged domain. In our case, we chose : $\sigma_c = 10$ MPa, $G_c=0.005$ mm.MPa and $p=0.1$.

\begin{figure}[ht!]
	\centering
	\subfloat[]{\includegraphics[width=0.47\textwidth]{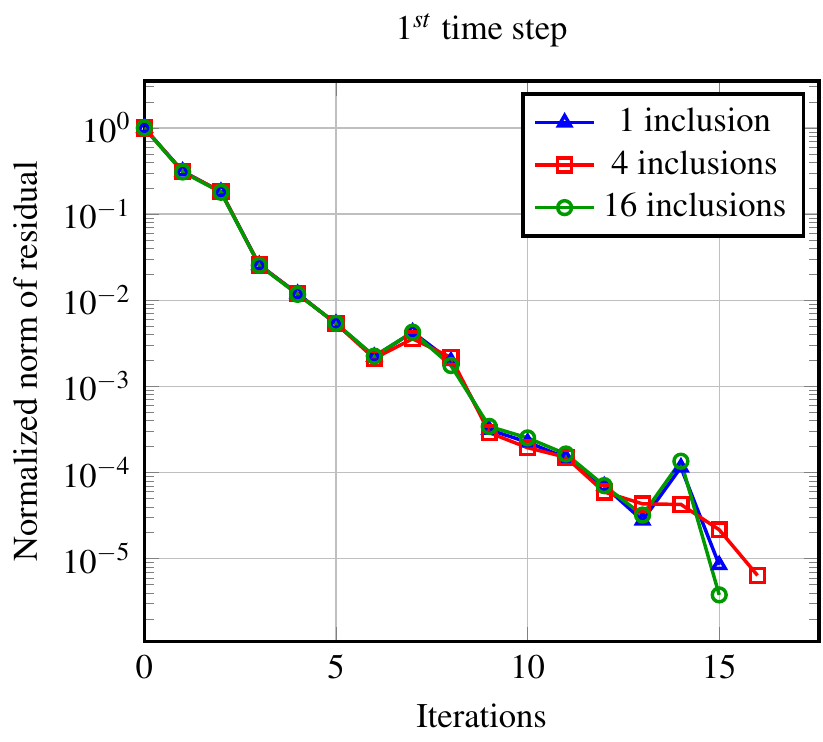}} \quad 
	\subfloat[]{\includegraphics[width=0.47\textwidth]{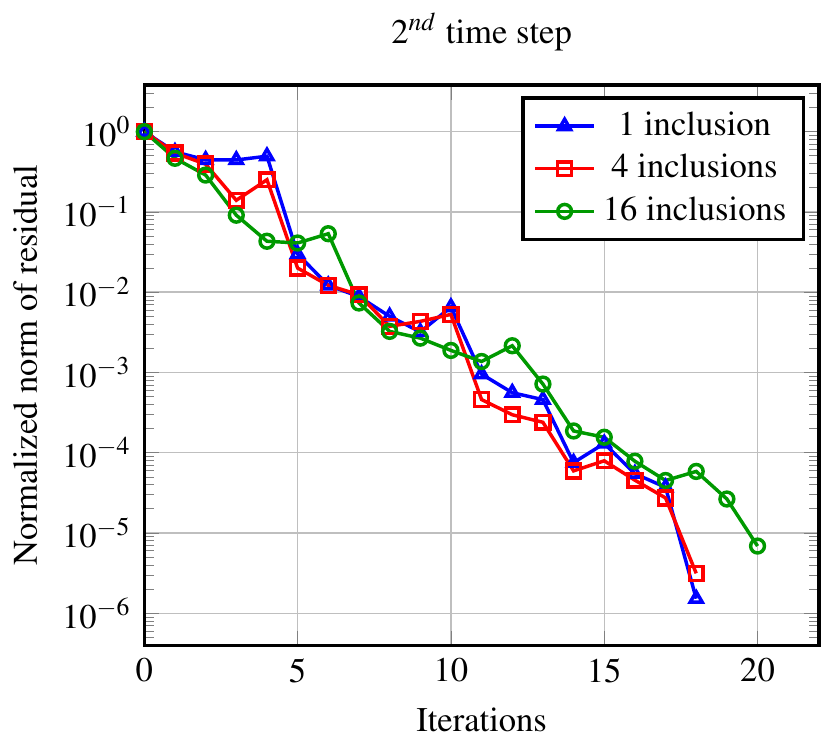}}\\
	\subfloat[]{\includegraphics[width=0.47\textwidth]{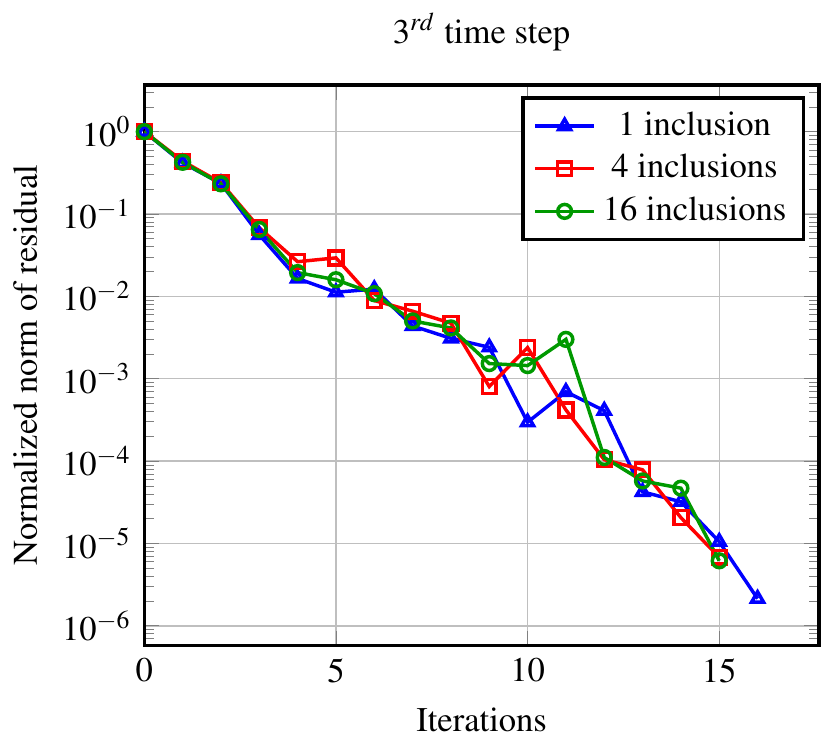}} \quad
	\subfloat[]{\includegraphics[width=0.47\textwidth]{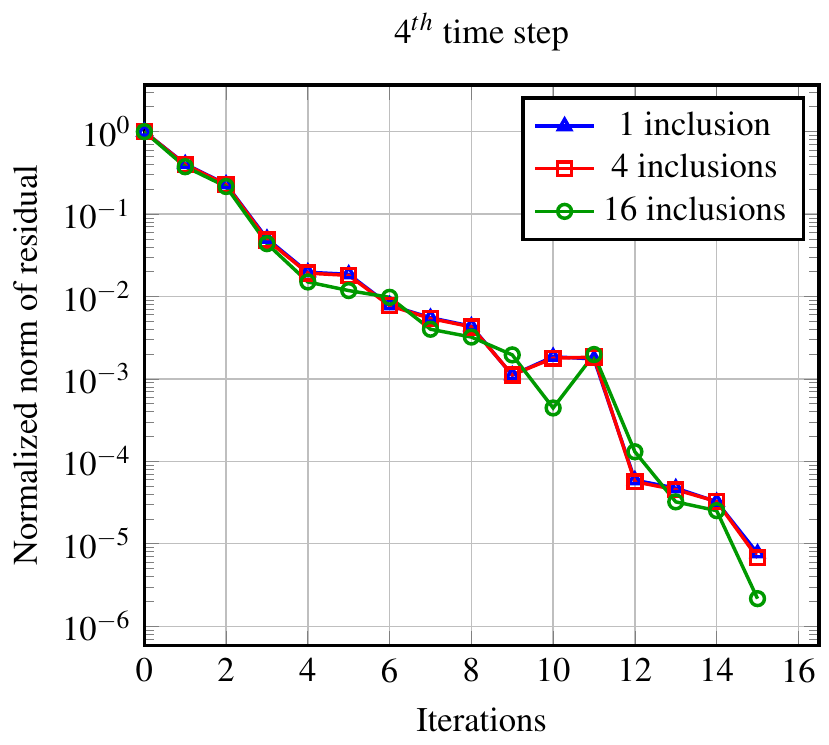}}
	\caption{Convergence of the non-invasive global/local algorithm for the plate problem with multiple inclusions and delamination. One graph corresponds to one loading step.}
	\label{fig:convInclusions}
\end{figure}

For the non-linear simulation, we considered 4 time steps with linear load increments to reach a final displacement of $W_d =0.0175$ mm on the right edge of the global structure. The problem was solved for each loading step using our non-invasive algorithm~\eqref{eq:nonint_global_dis_V3}-\eqref{eq:nonint_local_dis_V3} accelerated with the Aitken dynamic relaxation. In terms of results, we first show in Fig.~\ref{fig:convInclusions} that the convergence of the non-invasive global/local algorithm does not depend on the number of inclusions. Underlining that the different non-linear local problems can be naturally solved in parallel here, this accounts for the scalability of the algorithm. In other terms, our approach can be used as an efficient non-linear domain decomposition solver, as proposed in~\cite{Duval16} in standard FEM. The global model actually plays the role of the coarse problem of domain decomposition approaches, which allows to transmit directly the information all over the sub-domains. Furthermore, it may be noticed that the convergence of the non-invasive algorithm is similar for each loading step: between 10 and 15 iterations are sufficient to obtain a residual of $10^{-4}$. 

\begin{figure}
	\centering
	\subfloat[Stress $\sigma_{xx}$.]{\includegraphics[width=0.55\textwidth]{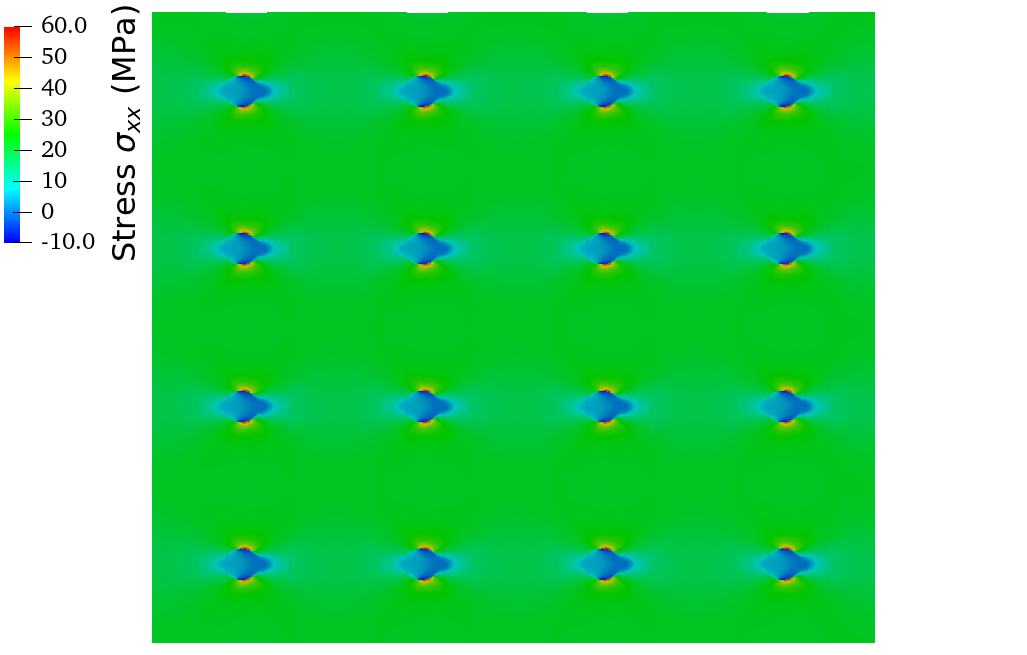}}
	\subfloat[Stress $\sigma_{yy}$.]{\includegraphics[width=0.55\textwidth]{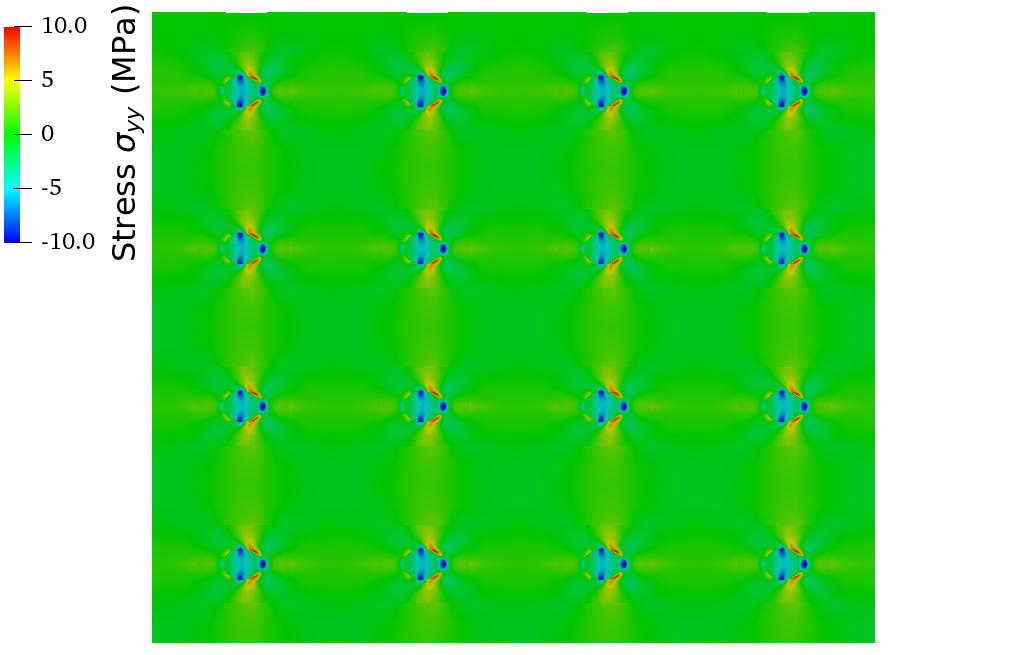}}
	\caption{Stress distribution  for the plate problem with 16 inclusions at the last loading step.}
	\label{fig:stresses_inclusions}
\end{figure}

\begin{figure}
	\centering
	\subfloat[Stress $\sigma_{yy}$.\label{fig:zoomStress}]{\includegraphics[width=0.55\textwidth]{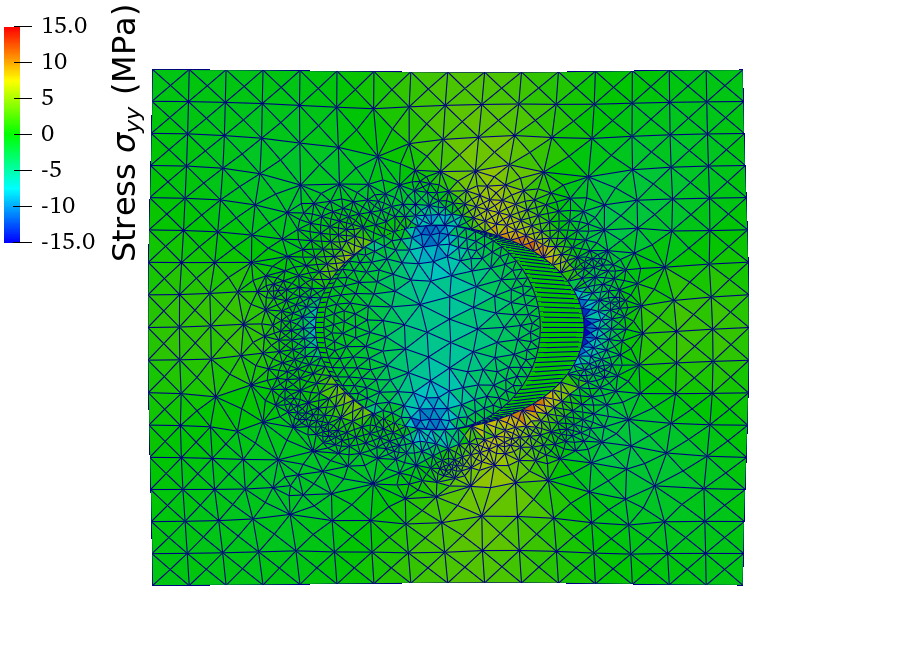}}
	\subfloat[Status of the damage in the cohesive elements - red : broken element, white : damaged element, blue : undamaged element.\label{fig:damage}]{\includegraphics[width=0.55\textwidth]{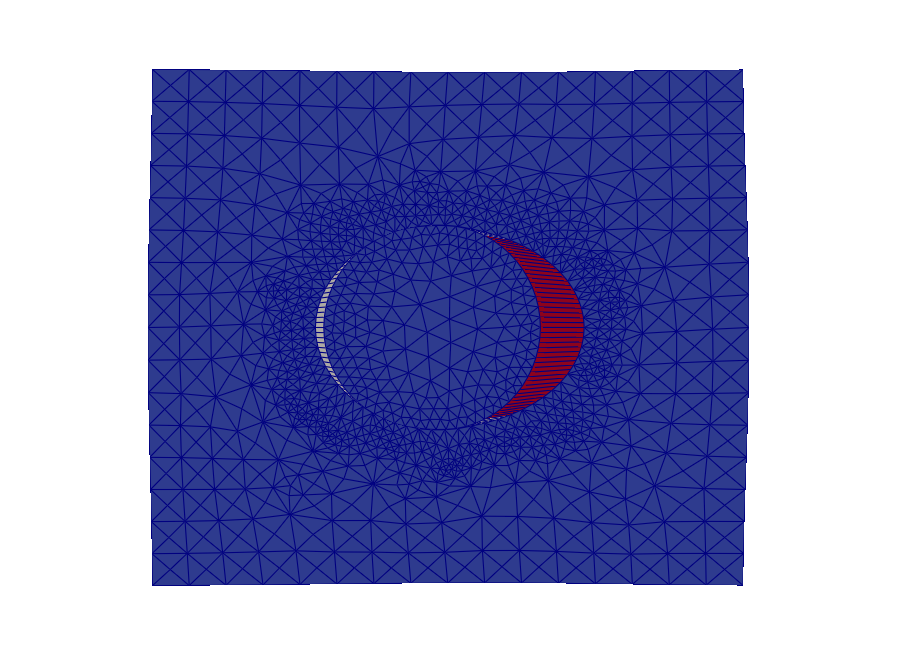}}
	\caption{Zoom on a local model for the plate problem with 16 inclusions at the last loading step.}
	\label{fig:zoom_damage}
\end{figure}

\begin{figure}
	\centering
	\includegraphics[width=0.6\textwidth]{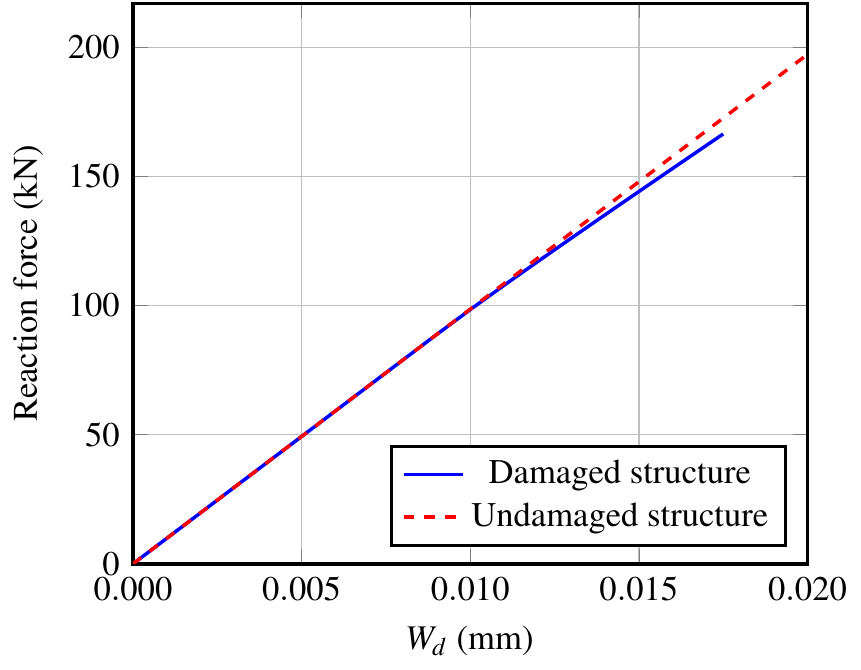}
	\caption{Reaction force versus prescribed displacement at each loading step for the plate problem with 16 inclusions.}
	\label{fig:reacInclusion}
\end{figure}

Then, Fig.~\ref{fig:stresses_inclusions} shows the stress distribution for the 16-inclusion problem at the last loading step. Stress concentrations can be observed around the inclusions which means that parts of the interfaces were severely damaged. This is further confirmed in Fig.~\ref{fig:zoom_damage} where a zoom over one inclusion is performed. The stress distribution can be better appreciated in Fig.~\ref{fig:zoomStress} and in Fig.~\ref{fig:damage} we plot the status of the damage in the cohesive elements. We see that almost all the interface were damaged and moreover the right part was almost completely broken. Finally, in Fig.~\ref{fig:reacInclusion}, we plot the reaction forces according to the prescribed displacement. The impact of the damage on the cohesive zones is noticed as the reaction force (blue curve) decreases compared to the curve with undamaged cohesive zones (dashed red curve). Let us note here that the load steps were quite large and this avoided numerical difficulties to capture some snap-back phenomena. Indeed with potential brutal loss of stiffness due to quick rupture of the cohesive zone, the precedent curve could involve snap-backs. With a finer load discretization, the Newton solver could fail and continuation techniques should be used. 

\subsection{3D mechanical assembly example with preload and frictional contact}

\begin{figure}
	\centering
	\includegraphics[width=0.7\textwidth]{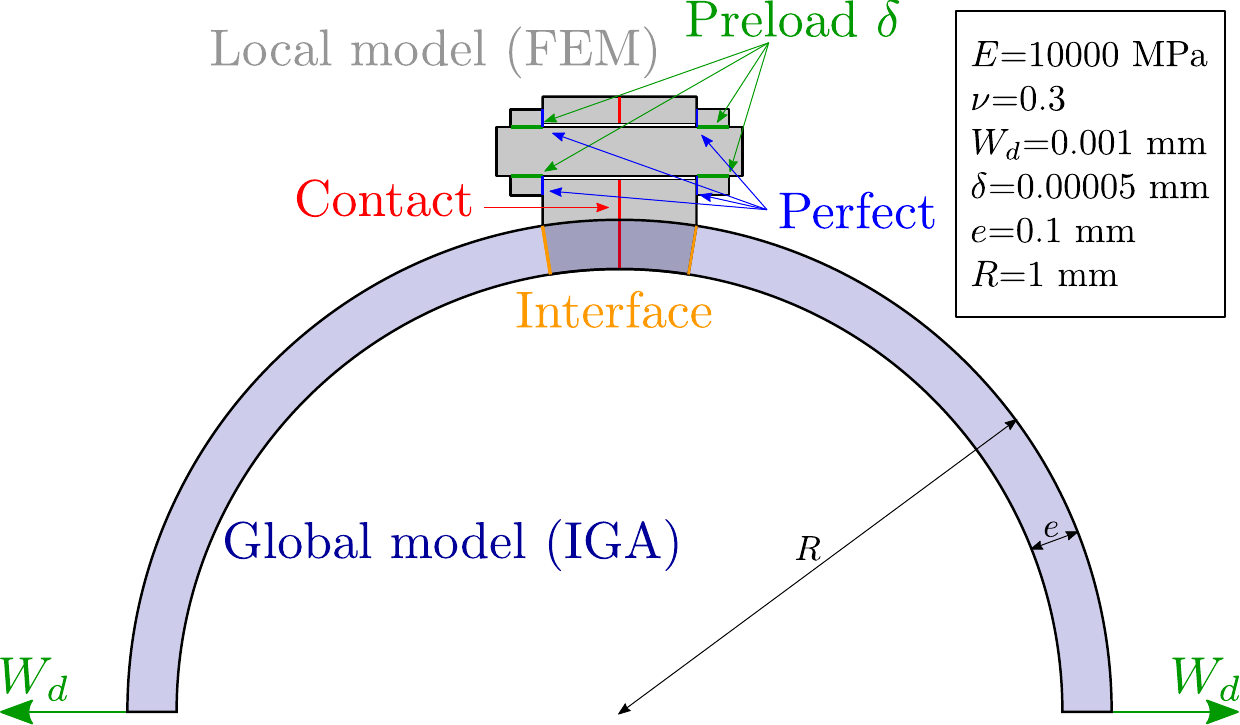}
	\caption{3D mechanical assembly example: description and data of the problem.}
	\label{fig:pb_assembly}
\end{figure}

\begin{figure}
	\centering
	\subfloat[Overview and zoomed window around the global/local interface. \label{fig:assembly_mesh_complete}]{\includegraphics[width=0.45\textwidth]{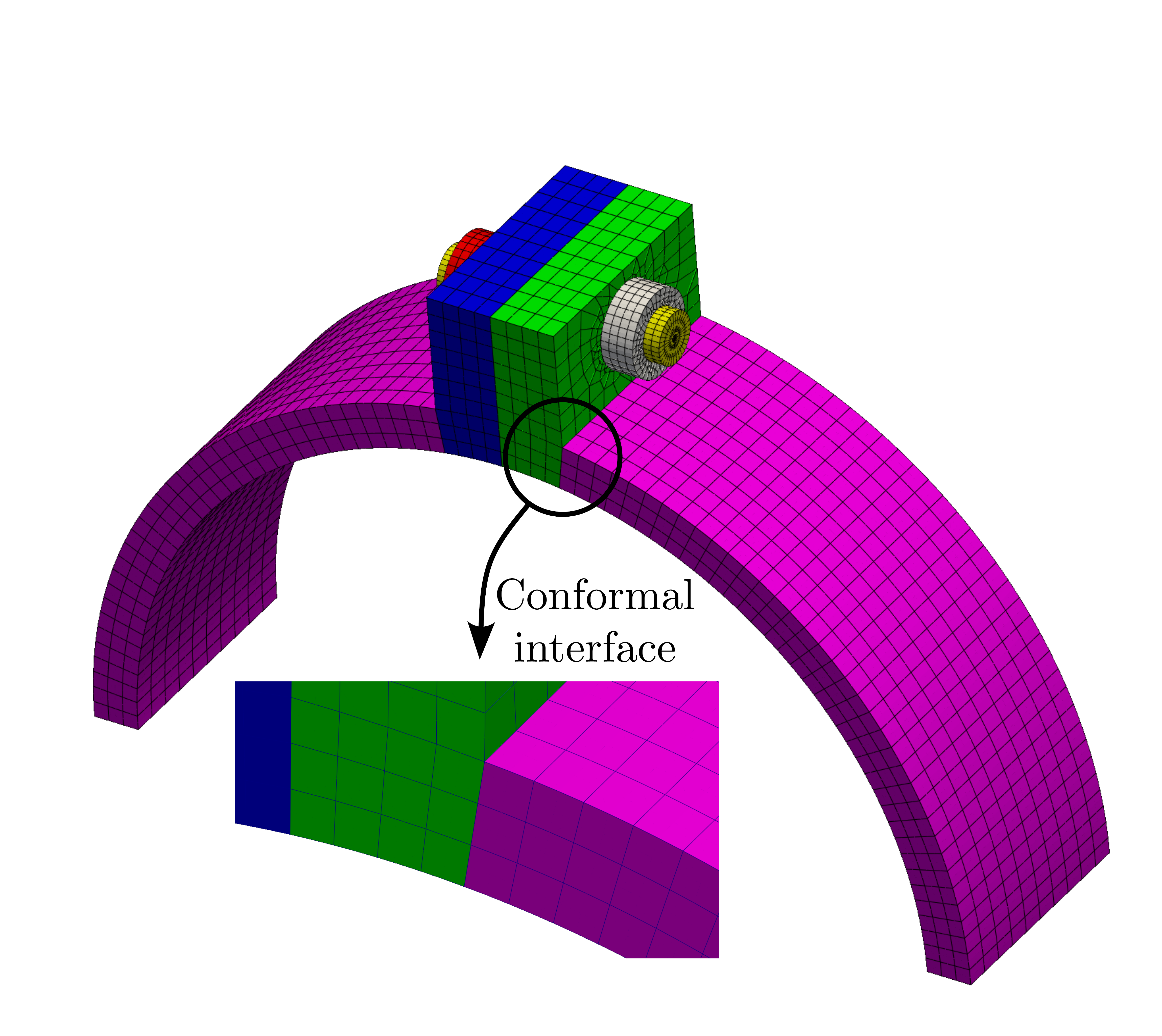}} \quad
	\subfloat[Discretization of the contact zone (a cut is performed in the between the two perforated plates in contact). \label{fig:mesh_contact_zone}]{\includegraphics[width=0.45\textwidth]{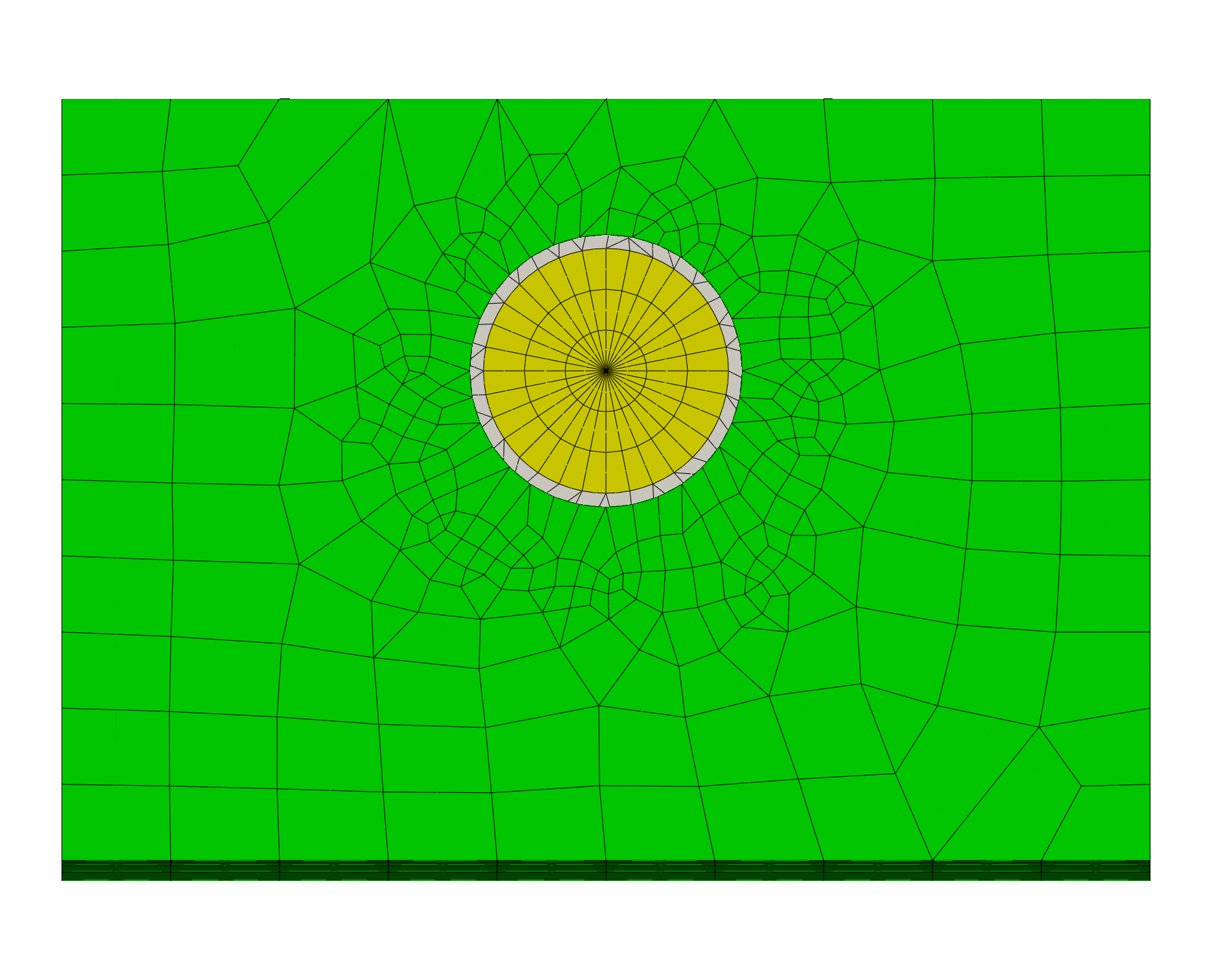}}
	\caption{Meshes for the 3D mechanical assembly example.}
	\label{fig:pb_assembly_mesh}
\end{figure}

As a last illustration, a 3D frictional contact problem of a bolted assembly was investigated. The test case along with all the model parameters are given in Fig.~\ref{fig:pb_assembly}. The purpose of this example is to show that the developed non-invasive hybrid procedure can be easily applied to solve general 3D contact problems. The assembly was composed of two quarter cylinders and one bolt that acts on two perforated plates which extend each of the pieces of cylinders. For the modeling, the global model was only made of a continuous half cylinder with quadratic NURBS elements (see blue region in Fig.~\ref{fig:pb_assembly}). This half cylinder was cut and bolted through a complex local FE model (see gray region in Fig.~\ref{fig:pb_assembly}). More precisely, the local model consisted in five solids: two perforated plates linked to the global model across the global/local interface, and a screw and two nuts that enforced the contact between the two plates thanks to the application of a preload in a pre-processing step. The dark blue region in Fig.~\ref{fig:pb_assembly} thus concerns domain $\Omega_{12}$, the NURBS solution is replaced by the local FE solution in this part. Overall, the local model included 5 interfaces (see again Fig.~\ref{fig:pb_assembly}): a frictional contact interface between the two plates (see red line), two preload interfaces between the nuts and the screw that were used to impose a relative displacement in the bolt (see green lines), and two perfect interfaces between the plates and the nuts (see blue lines). In the end, these five interfaces produced tension in the screw and compression in the plates. As for the discretization, quadratic tetrahedrons and pyramids were used for the local model and the strategy depicted in Fig.~\ref{fig:gIGAlFEM_papier_principe} was applied in 3D to recover a conforming global/local interface (see Fig.~\ref{fig:pb_assembly_mesh}). Eventually, displacement boundary conditions were prescribed at the bottom of the global model (half cylinder) to extend the bolted assembly.

\begin{figure}
	\centering
	\includegraphics[scale=0.9]{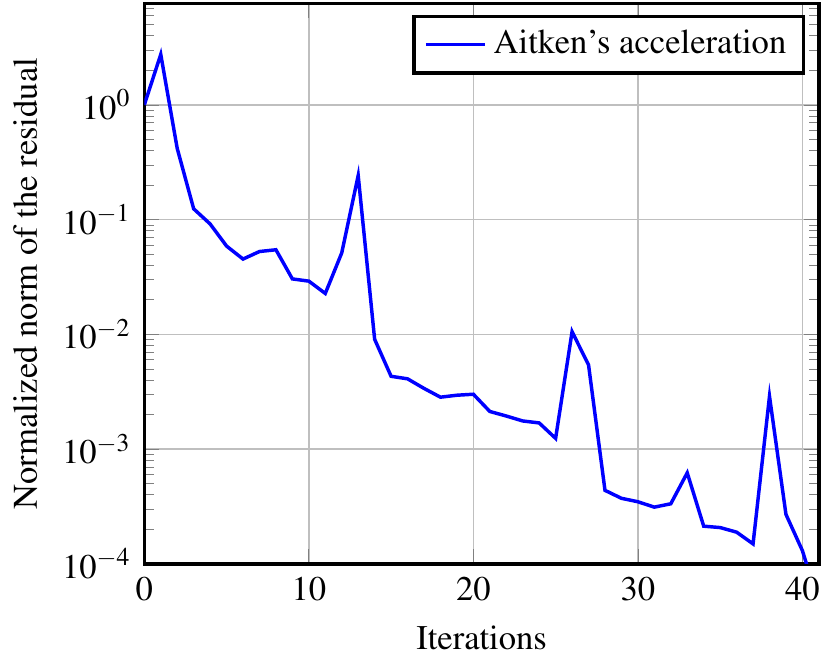}
	\caption{Convergence of the non-invasive algorithm for the 3D mechanical assembly test case. }
	\label{fig:convBoulon}
\end{figure}

\begin{figure}
	\centering
	\includegraphics[width=0.9\textwidth]{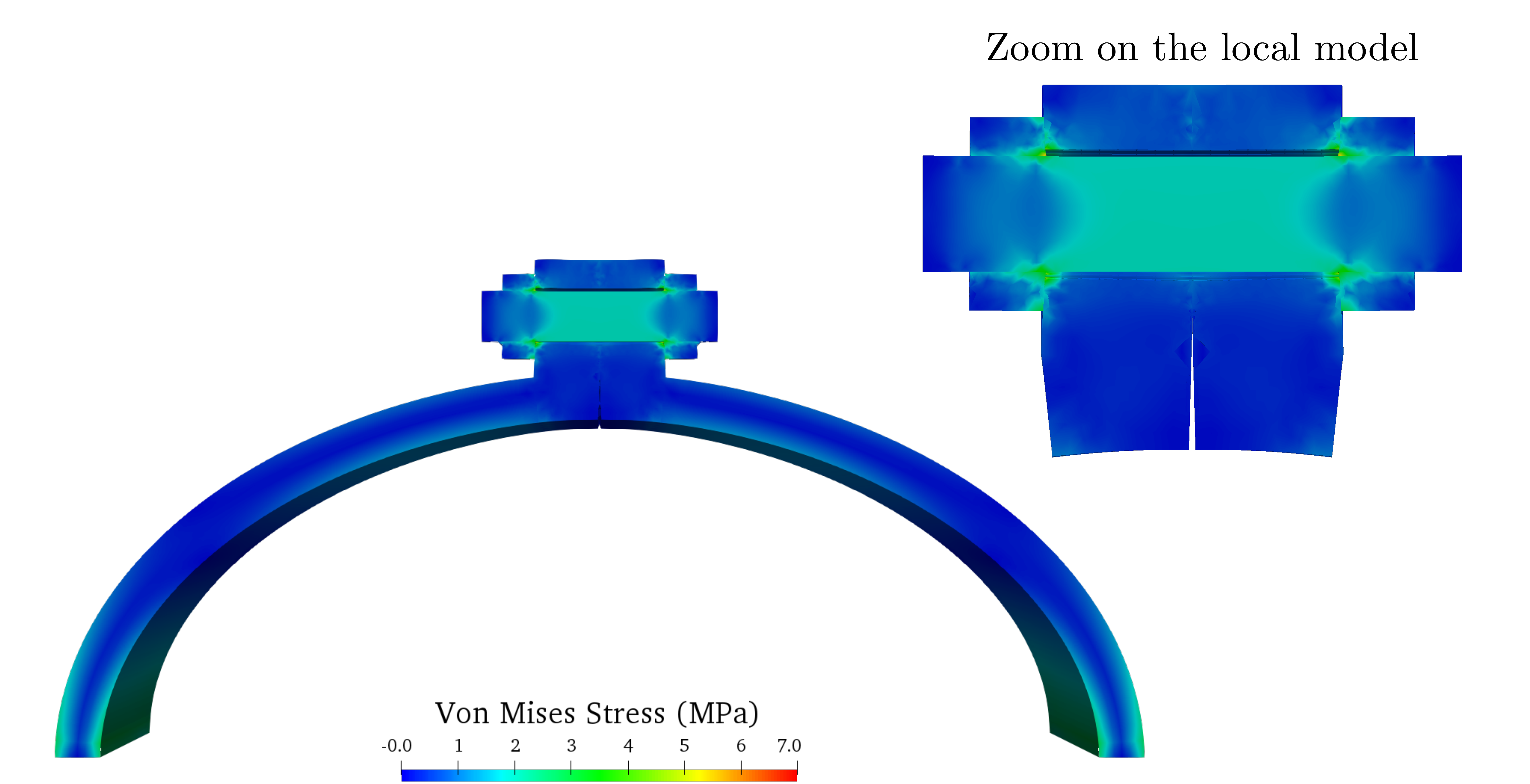}
	\caption{Obtained Von Mises stress and deformed configuration (scale factor 200). Half of the structure is shown (a cut is performed in the symmetry plan of the structure).}
	\label{fig:vonMisesBoulon}
\end{figure}

\begin{figure}
	\centering
	\includegraphics[width=0.7\textwidth]{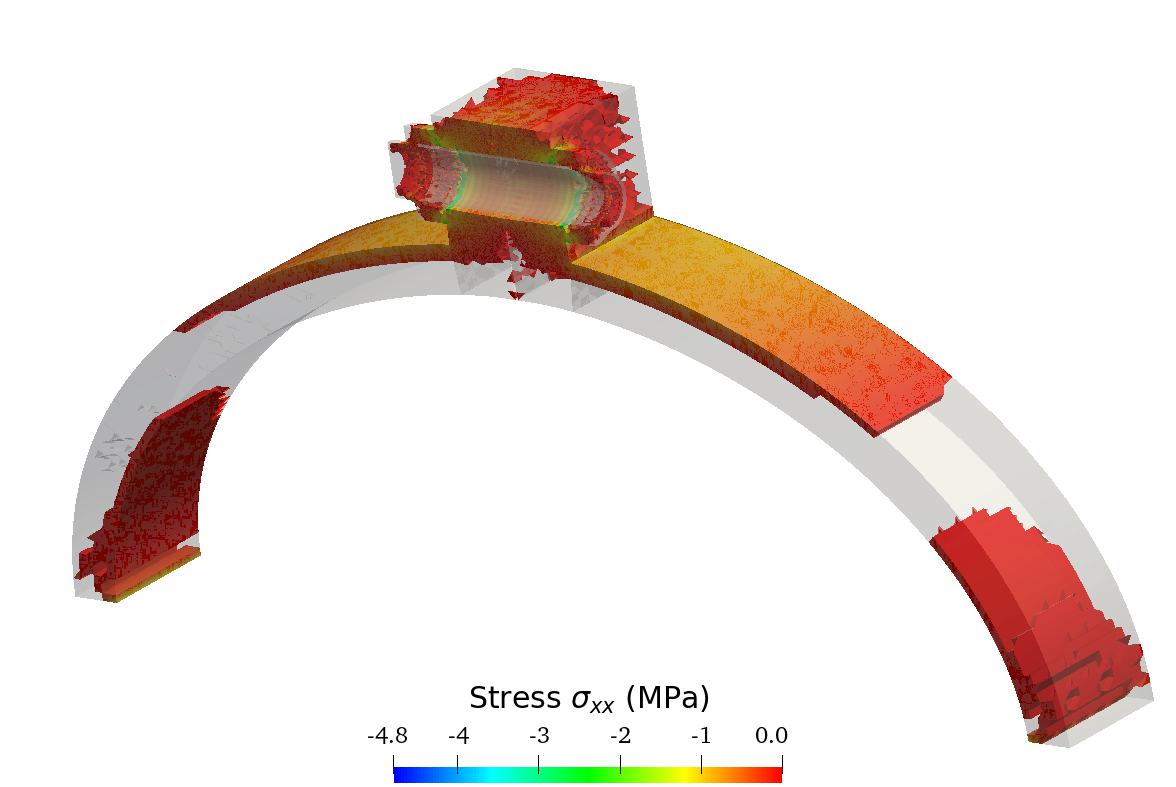}
	\caption{Highlighting of the compression zones (\emph{i.e.} where $\sigma_{xx}$ is negative).}
	\label{fig:sXXboulon}
\end{figure}

The convergence of the non-invasive algorithm is provided in Fig~\ref{fig:convBoulon}. An Aitken's acceleration was necessary to obtain a residual of $10^{-4}$ in a few tens of iterations as in the previous non-linear test cases. Then, the deformed shape along with the Von Mises stress is depicted in Fig.~\ref{fig:vonMisesBoulon}. Once again, a smooth transition of the stress is observed at the global/local interface which confirms the performance of our non-invasive coupling scheme in 3D. Moreover, due to the applied Dirichlet boundary conditions, a slight detachment can be observed between the two quarter cylinders near the interior radius whereas the bolt maintains the contact around the screw. Obviously, stress concentrations can also be reported close to the nuts. Finally, we show in Fig.~\ref{fig:sXXboulon} the parts of the assembly that are in compression (\emph{i.e.} where $\sigma_{xx}$ is negative). As expected, the screw is in tension due to the preload and leads to a large compression zone under the nuts, which helps keeping the plates in contact. This good behavior of this last case demonstrates the potential of our method to treat more representative applications.

\begin{remark}
	Let us notice that several bolts could be easily considered using exactly the same strategy but with several local models, in the previous test case. The superiority of the non-invasive approach would be even clearer from a modeling and a solution point of view. Indeed, solving a multi-contact problem in a monolithic way is not an easy task as convergence properties of the non-linear solvers are worsened by the increase of the number of contact surfaces. Here, the different contact problems associated with each local models could be naturally solved in parallel, therefore we would end up with a non-invasive, non-linear Domain Decomposition strategy (see again~\cite{Duval16} for more information on this topic).
\end{remark}

\section{Conclusion}
\label{sec:conclu}

In this work, we managed to implement the coupling between IGA and FEM in a seamless and automatic manner for the robust and efficient multiscale global/local simulation of structures. IGA is used for the global model since it allows for a proper description of the global geometry and for an increased per-DOF accuracy to capture global, regular responses. Conversely, FEM is considered for the simulation of the local behavior since its meshing capabilities allow for the proper description of geometrical details (such as cracks) and its low regular discretization space appears more adapted to the modeling of local, strongly non-linear or even singular phenomena (\emph{e.g.}, contact, delamination, heterogeneities). From this point of view, our strategy can be viewed as an improvement of a full global/local FEM in the sense that it fosters a direct link with CAD for the global geometry and it saves many DOF for the same accuracy in the global region. Reciprocally, our method also enhances global/local IGA from the perspective that it facilitates the modelling of complex local behaviors within an IG (possible multi) patch structure by making use of existing, optimized FE routines specifically developed for this purpose.

From a technical viewpoint, our strategy  falls in the class of so-called non-invasive global/local algorithms that is currently gathering an important momentum in the community. The method relies on an iterative Dirichlet-to-Neumann process where the initial global IG model to be enriched is never modified. This constitutes the first key point in IGA to reach a non-invasive strategy since it avoids costly spline re-parametrization procedures that may have been necessary otherwise to incorporate a free local region. In addition, it has to be noted that since only the initial global IG operator is involved, the IG system to be solved remains well conditioned regardless of the shape of the local region. Then, to achieve full non-invasiveness, we proposed to resort to the existing FEM-to-IGA bridge, based on Bézier or Lagrange extraction, to transform the interface within the initial global IG model into a FE interface on which the local FE mesh can be constructed by calling upon efficient (classic) meshing procedures. The situation of a conforming coupling interface was thus reached and we were able to express the hybrid IGA/FEM coupling by means of only the FEM-to-IGA operator (easily built from any spline libraries) and standard FE trace operators that are available in (possibly industrial) FE codes. It results that any robust FE code suitable for the modelling of complex local behaviors can be used in a plug-and-play manner. Finally, depending on how deeply the FEM-to-IGA operator is used, our implementation is generic in terms of programming environments: the users may have in hand an IG code (performing standard elasticity) and wish to couple it with a specific FE software to model complex local phenomena, or the users only have at their disposal FE packages.

For demonstration purpose, we considered the second situation for the numerical experiments which is the most extreme one. The implementation was performed using the open-source FE industrial software package Code\_Aster~\cite{aster} developed by the EDF R\&D company. More precisely, the strategy was first validated through a simple 2D linear elastic problem where an analytical solution was available, which allowed to properly account for the accuracy of the non-invasive coupling. Then, we highlighted the flexibility and robustness of the methodology by considering the enrichment of a global 2D IG linear elastic model by means of (i) a specific FE mesh incorporating holes, cracks and frictional contact between the different lips of the cracks, and (ii) several FE local meshes including inclusions along with delamination on their boundary in view of simulating fiber-reinforced composites. The latter test case also enabled to show another attractive property of the non-invasive algorithm: it naturally leads to a scalable non-linear domain decomposition solver when several local models are incorporated. Finally, a 3D mechanical assembly problem was solved where the local FE mesh was used to model the bolt while the IG model represented a half cylinder.

Although restricted in this contribution to standard B-spline and NURBS technologies, it may be underlined that our method easily extends to recent advanced splines since most of them support Bézier extraction (\emph{e.g.}, hierarchical B-Splines and NURBS~\cite{Hennig16,Angella20}, and hierarchical T-Splines~\cite{Evans15,Chen18}). Such an extension would be interesting in case of a tiny local model within a large global structure. Rather than refining the global IG model everywhere, advanced splines that allow for local mesh refinement could be used in order to forward more efficiently information from the local zone to the rest of the structure. Then, it is obvious that if one has an optimized IG code in hand, it is more efficient to directly use it than recovering IGA from a FE code, as done here with Code\_Aster, since it removes the projection step. In this sense, our method is consistent with the current development of open-source IG libraries  (see, \emph{e.g.}, Nutils (http://www.nutils.org) or pyiga (https://github.com/c-f-h/pyiga) as Python libraries, tIGAr~\cite{Kamensky19} for a FENICS-based implementation, and YETI (https://lamcosplm.insa-lyon.fr/projects/yeti/) for a Fortran-Python implementation). Even more interestingly, since the IG code is only used at the global scale, where standard elasticity may be sufficient, our strategy also appears totally relevant with the newly developed fast assembly and solution procedures for IGA in the linear regime (see, \emph{e.g.}, sum factorization~\cite{Antolin15}, use of look-up tables~\cite{Mantzaflaris15,Hirschler22}, weighted quadrature~\cite{Hiemstra19}, and domain decomposition solvers~\cite{Hirschler19,Bosy20}).

\end{document}